

\ifx\begin\undefined\else\global\advance\srcdepth by
1\expandafter \fi

\def\begin{}
\newcount\srcdepth
\srcdepth=1

\outer\def\bye{\global\advance\srcdepth by -1
  \ifnum\srcdepth=0
    \def\endcmd{\vfill\eject\nopagenumbers\par\vfill\supereject\end}
  \else\def\endcmd{}\fi
  \endcmd
}


\baselineskip=13pt
\hsize = 5.5truein
\hoffset = 0.5truein
\vsize = 8.5truein
\voffset = 0.2truein
\emergencystretch = 0.05\hsize

\overfullrule=0pt

\newif\ifblackboardbold

\blackboardboldtrue


\font\sectionfont=cmbx12


\newfam\bboldfam
\ifblackboardbold
\font\tenbbold=msbm10
\font\sevenbbold=msbm7
\font\fivebbold=msbm5
\textfont\bboldfam=\tenbbold
\scriptfont\bboldfam=\sevenbbold
\scriptscriptfont\bboldfam=\fivebbold
\def\bbold{\fam\bboldfam\tenbbold}
\else
\def\bbold{\bf}
\fi


\font\Arm=cmr8
\font\Ai=cmmi8
\font\Asy=cmsy8
\font\Abf=cmbx8
\font\Brm=cmr6
\font\Bi=cmmi6
\font\Bsy=cmsy6
\font\Bbf=cmbx6
\font\Crm=cmr5
\font\Ci=cmmi5
\font\Csy=cmsy5
\font\Cbf=cmbx5

\ifblackboardbold
\font\Abbold=msbm10 at 8pt
\font\Bbbold=msbm7 at 6pt
\font\Cbbold=msbm5
\fi

\def\smallmath{%
\textfont0=\Arm \scriptfont0=\Brm \scriptscriptfont0=\Crm
\textfont1=\Ai \scriptfont1=\Bi \scriptscriptfont1=\Ci
\textfont2=\Asy \scriptfont2=\Bsy \scriptscriptfont2=\Csy
\textfont\bffam=\Abf \scriptfont\bffam=\Bbf \scriptscriptfont\bffam=\Cbf
\def\rm{\fam0\Arm}\def\mit{\fam1}\def\oldstyle{\fam1\Ai}%
\def\bf{\fam\bffam\Abf}%
\ifblackboardbold
\textfont\bboldfam=\Abbold
\scriptfont\bboldfam=\Bbbold
\scriptscriptfont\bboldfam=\Cbbold
\def\bbold{\fam\bboldfam\Abbold}%
\fi
}








\newlinechar=`@
\def\forwardmsg#1#2#3{\immediate\write16{@*!*!*!* forward reference should
be: @\noexpand\forward{#1}{#2}{#3}@}}
\def\nodefmsg#1{\immediate\write16{@*!*!*!* #1 is an undefined reference@}}

\def\forwardsub#1#2{\def\newref{{#2}{#1}}}

\def\forward#1#2#3{%
\expandafter\expandafter\expandafter\forwardsub\expandafter{#3}{#2}
\expandafter\ifx\csname#1\endcsname\relax\else%
\expandafter\ifx\csname#1\endcsname\newref\else%
\forwardmsg{#1}{#2}{#3}\fi\fi%
\expandafter\let\csname#1\endcsname\newref}

\def\firstarg#1{\expandafter\argone #1}\def\argone#1#2{#1}
\def\secondarg#1{\expandafter\argtwo #1}\def\argtwo#1#2{#2}

\def\ref#1{\expandafter\ifx\csname#1\endcsname\relax
  {\nodefmsg{#1}\bf`#1'}\else
  \expandafter\firstarg\csname#1\endcsname
  ~\expandafter\secondarg\csname#1\endcsname\fi}

\def\refs#1{\expandafter\ifx\csname#1\endcsname\relax
  {\nodefmsg{#1}\bf`#1'}\else
  \expandafter\firstarg\csname #1\endcsname
  s~\expandafter\secondarg\csname#1\endcsname\fi}

\def\refn#1{\expandafter\ifx\csname#1\endcsname\relax
  {\nodefmsg{#1}\bf`#1'}\else
  \expandafter\secondarg\csname #1\endcsname\fi}



\def\widow#1{\vskip 0pt plus#1\vsize\goodbreak\vskip 0pt plus-#1\vsize}



\def\marginlabel#1{}

\def\showlabelsabove{
\font\labelfont=cmss10 at 6pt
\def\marginlabel##1{\rlap{\smash{\raise 10pt\hbox{\labelfont##1}}}}
}

\newcount\seccount
\newcount\proccount
\seccount=0
\proccount=0

\def\stdskip{\vskip 9pt plus3pt minus 3pt}
\def\stdbreak{\par\removelastskip\penalty-100\stdskip}

\def\proof{\stdbreak\noindent{\sl Proof. }}

\def\qed{\vrule height 1.2ex width .9ex depth .1ex}

\def\Box{
  \ifmmode\eqno\qed
  \else\ifvmode\removelastskip\line{\hfil\qed}
  \else\unskip\quad\hskip-\hsize
    \hbox{}\hskip\hsize minus 1em\qed\par
  \fi\stdbreak\fi}

\def\references{
  \removelastskip
  \widow{.05}
  \vskip 24pt plus 6pt minus 6 pt
  \leftline{\sectionfont References}
  \nobreak\stdskip\noindent}

\def\ifempty#1#2\endB{\ifx#1\endA}
\def\makeref#1#2#3{\ifempty#1\endA\endB\else\forward{#1}{#2}{#3}\fi}

\outer\def\section#1 #2\par{
  \removelastskip
  \global\advance\seccount by 1
  \global\proccount=0\relax
                \edef\numtoks{\number\seccount}
  \makeref{#1}{Section}{\numtoks}
  \widow{.05}
  \vskip 24pt plus 6pt minus 6 pt
  \message{#2}
  \leftline{\marginlabel{#1}\sectionfont\numtoks\quad #2}
  \nobreak\stdskip}

\def\proclamation#1#2{
  \outer\expandafter\def\csname#1\endcsname##1 ##2\par{
  \stdbreak
  \advance\proccount by 1
  \edef\numtoks{\number\seccount.\number\proccount}
  \makeref{##1}{#2}{\numtoks}
  \noindent{\marginlabel{##1}\bf #2 \numtoks\enspace}
  {\sl##2\par}
  \stdbreak}}

\def\othernumbered#1#2{
  \outer\expandafter\def\csname#1\endcsname##1{
  \stdbreak
  \advance\proccount by 1
  \edef\numtoks{\number\seccount.\number\proccount}
  \makeref{##1}{#2}{\numtoks}
  \noindent{\marginlabel{##1}\bf #2 \numtoks\enspace}}}

\proclamation{definition}{Definition}
\proclamation{lemma}{Lemma}
\proclamation{proposition}{Proposition}
\proclamation{theorem}{Theorem}
\proclamation{corollary}{Corollary}
\proclamation{conjecture}{Conjecture}

\othernumbered{example}{Example}
\othernumbered{remark}{Remark}
\othernumbered{construction}{Construction}






\newcount\figcount
\figcount=0
\newbox\drawing
\newcount\drawbp
\newdimen\drawx
\newdimen\drawy
\newdimen\ngap
\newdimen\sgap
\newdimen\wgap
\newdimen\egap

\def\drawbox#1#2#3{\vbox{
  \setbox\drawing=\vbox{\offinterlineskip\epsfbox{#2.eps}\kern 0pt}
  \drawbp=\epsfurx
  \advance\drawbp by-\epsfllx\relax
  \multiply\drawbp by #1
  \divide\drawbp by 100
  \drawx=\drawbp truebp
  \ifdim\drawx>\hsize\drawx=\hsize\fi
  \epsfxsize=\drawx
  \setbox\drawing=\vbox{\offinterlineskip\epsfbox{#2.eps}\kern 0pt}
  \drawx=\wd\drawing
  \drawy=\ht\drawing
  \ngap=0pt \sgap=0pt \wgap=0pt \egap=0pt 
  \setbox0=\vbox{\offinterlineskip
    \box\drawing \ifgridlines\drawgrid\drawx\drawy\fi #3}
  \kern\ngap\hbox{\kern\wgap\box0\kern\egap}\kern\sgap}}

\def\draw#1#2#3{
  \setbox\drawing=\drawbox{#1}{#2}{#3}
  \advance\figcount by 1
  \goodbreak
  \midinsert
  \centerline{\ifgridlines\boxgrid\drawing\fi\box\drawing}
  \smallskip
  \vbox{\offinterlineskip
    \centerline{Figure~\number\figcount}
    \smash{\marginlabel{#2}}}
  \endinsert}

\def\nextfigtoks{%
  \advance\figcount by 1%
  \edef\numtoks{\number\figcount}%
  \advance\figcount by -1}

\newif\ifgridlines
\newbox\figtbox
\newbox\figgbox
\newdimen\figtx
\newdimen\figty

\newdimen\bwd
\bwd=2sp 

\def\hline#1{\vbox{\smash{\hbox to #1{\leaders\hrule height \bwd\hfil}}}}

\def\vline#1{\hbox to 0pt{%
  \hss\vbox to #1{\leaders\vrule width \bwd\vfil}\hss}}

\def\clap#1{\hbox to 0pt{\hss#1\hss}}
\def\vclap#1{\vbox to 0pt{\offinterlineskip\vss#1\vss}}

\def\hstutter#1#2{\hbox{%
  \setbox0=\hbox{#1}%
  \hbox to #2\wd0{\leaders\box0\hfil}}}

\def\vstutter#1#2{\vbox{
  \setbox0=\vbox{\offinterlineskip #1}
  \dp0=0pt
  \vbox to #2\ht0{\leaders\box0\vfil}}}

\def\crosshairs#1#2{
  \dimen1=.002\drawx
  \dimen2=.002\drawy
  \ifdim\dimen1<\dimen2\dimen3\dimen1\else\dimen3\dimen2\fi
  \setbox1=\vclap{\vline{2\dimen3}}
  \setbox2=\clap{\hline{2\dimen3}}
  \setbox3=\hstutter{\kern\dimen1\box1}{4}
  \setbox4=\vstutter{\kern\dimen2\box2}{4}
  \setbox1=\vclap{\vline{4\dimen3}}
  \setbox2=\clap{\hline{4\dimen3}}
  \setbox5=\clap{\copy1\hstutter{\box3\kern\dimen1\box1}{6}}
  \setbox6=\vclap{\copy2\vstutter{\box4\kern\dimen2\box2}{6}}
  \setbox1=\vbox{\offinterlineskip\box5\box6}
  \smash{\vbox to #2{\hbox to #1{\hss\box1}\vss}}}

\def\boxgrid#1{\rlap{\vbox{\offinterlineskip
  \setbox0=\hline{\wd#1}
  \setbox1=\vline{\ht#1}
  \smash{\vbox to \ht#1{\offinterlineskip\copy0\vfil\box0}}
  \smash{\vbox{\hbox to \wd#1{\copy1\hfil\box1}}}}}}

\def\drawgrid#1#2{\vbox{\offinterlineskip
  \dimen0=\drawx
  \dimen1=\drawy
  \divide\dimen0 by 10
  \divide\dimen1 by 10
  \setbox0=\hline\drawx
  \setbox1=\vline\drawy
  \smash{\vbox{\offinterlineskip
    \copy0\vstutter{\kern\dimen1\box0}{10}}}
  \smash{\hbox{\copy1\hstutter{\kern\dimen0\box1}{10}}}}}

\def\figtext#1#2#3#4#5{
  \setbox\figtbox=\hbox{#5}
  \dp\figtbox=0pt
  \figtx=-#3\wd\figtbox \figty=-#4\ht\figtbox
  \advance\figtx by #1\drawx \advance\figty by #2\drawy
  \dimen0=\figtx \advance\dimen0 by\wd\figtbox \advance\dimen0 by-\drawx
  \ifdim\dimen0>\egap\global\egap=\dimen0\fi
  \dimen0=\figty \advance\dimen0 by\ht\figtbox \advance\dimen0 by-\drawy
  \ifdim\dimen0>\ngap\global\ngap=\dimen0\fi
  \dimen0=-\figtx
  \ifdim\dimen0>\wgap\global\wgap=\dimen0\fi
  \dimen0=-\figty
  \ifdim\dimen0>\sgap\global\sgap=\dimen0\fi
  \smash{\rlap{\vbox{\offinterlineskip
    \hbox{\hbox to \figtx{}\ifgridlines\boxgrid\figtbox\fi\box\figtbox}
    \vbox to \figty{}
    \ifgridlines\crosshairs{#1\drawx}{#2\drawy}\fi
    \kern 0pt}}}}


\def\hpad#1#2#3{\hbox{\kern #1\hbox{#3}\kern #2}}
\def\vpad#1#2#3{\setbox0=\hbox{#3}\dp0=0pt\vbox{\kern #1\box0\kern #2}}



\def\stack#1#2#3{\vbox{\offinterlineskip
  \setbox2=\hbox{#2}
  \setbox3=\hbox{#3}
  \dimen0=\ifdim\wd2>\wd3\wd2\else\wd3\fi
  \hbox to \dimen0{\hss\box2\hss}
  \kern #1
  \hbox to \dimen0{\hss\box3\hss}}}


\def\hexp#1{%
  \setbox0=\hbox{${}^{#1}$}%
  \hbox to .5\wd0{\box0\hss}}



\def\bmatrix#1#2{{\smallmath\left[\vcenter{\halign
  {&\kern#1\hfil$##\mathstrut$\kern#1\cr#2}}\right]}}

\def\rightarrowmat#1#2#3{
  \setbox1=\hbox{\kern#2$\bmatrix{#1}{#3}$\kern#2}
  \,\vbox{\offinterlineskip\hbox to\wd1{\hfil\copy1\hfil}
    \kern 3pt\hbox to\wd1{\rightarrowfill}}\,}

\def\leftarrowmat#1#2#3{
  \setbox1=\hbox{\kern#2$\bmatrix{#1}{#3}$\kern#2}
  \,\vbox{\offinterlineskip\hbox to\wd1{\hfil\copy1\hfil}
    \kern 3pt\hbox to\wd1{\leftarrowfill}}\,}

\def\rightarrowbox#1#2{
  \setbox1=\hbox{\kern#1\hbox{\smallmath #2}\kern#1}
  \,\vbox{\offinterlineskip\hbox to\wd1{\hfil\copy1\hfil}
    \kern 3pt\hbox to\wd1{\rightarrowfill}}\,}

\def\leftarrowbox#1#2{
  \setbox1=\hbox{\kern#1\hbox{\smallmath #2}\kern#1}
  \,\vbox{\offinterlineskip\hbox to\wd1{\hfil\copy1\hfil}
    \kern 3pt\hbox to\wd1{\leftarrowfill}}\,}






\def\bookletdims{
  \hsize=5.25truein
  \vsize=7truein
}

\def\legalbooklet#1{
  \input quire
  \bookletdims
  \htotal=7.0truein
  \vtotal=8.5truein
  \hoffset=\htotal
  \advance\hoffset by -\hsize
  \divide\hoffset by 2
  \voffset=\vtotal
  \advance\voffset by -\vsize
  \divide\voffset by 2
  \advance\voffset by -.0625truein
  \shhtotal=2\htotal
  \horigin=0.0truein
  \vorigin=0.0truein
  \shstaplewidth=0.01pt
  \shstaplelength=0.66truein
  \shthickness=0pt
  \shoutline=0pt
  \shcrop=0pt
  \shvoffset=-1.0truein
  \ifnum#1>0\quire{#1}\else\qtwopages\fi
}

\def\preview{
  \input quire
  \bookletdims
  \hoffset=0.1truein
  \vtotal=8.5truein
  \shhtotal=14truein
  \voffset=\vtotal
  \advance\voffset by -\vsize
  \divide\voffset by 2
  \advance\voffset by -.0625truein
  \htotal=2\hoffset
  \advance\htotal by \hsize
  \horigin=0.0truein
  \vorigin=0.0truein
  \shstaplewidth=0.5pt
  \shstaplelength=0.5\vtotal
  \shthickness=0pt
  \shoutline=0pt
  \shcrop=0pt
  \shvoffset=-1.0truein
  \qtwopages
}

\def\twoup{
  \input quire
  \hsize=4.79452truein 
  \vsize=7truein
  \vtotal=8.5truein
  \shhtotal=11truein
  \hoffset=-2\hsize
  \advance\hoffset by \shhtotal
  \divide\hoffset by 6
  \voffset=\vtotal
  \advance\voffset by -\vsize
  \divide\voffset by 2
  \advance\voffset by -12truept
  \htotal=2\hoffset
  \advance\htotal by \hsize
  \horigin=0.0truein
  \vorigin=0.0truein
  \shstaplewidth=0.01pt
  \shstaplelength=0pt
  \shthickness=0pt
  \shoutline=0pt
  \shcrop=0pt
  \shvoffset=-1.0truein
  \qtwopages
}


\newcount\countA
\newcount\countB
\newcount\countC

\def\monthname{\begingroup
  \ifcase\number\month
    \or January\or February\or March\or April\or May\or June\or
    July\or August\or September\or October\or November\or December\fi
\endgroup}

\def\dayname{\begingroup
  \countA=\number\day
  \countB=\number\year
  \advance\countA by 0 
  \advance\countA by \ifcase\month\or
    0\or 31\or 59\or 90\or 120\or 151\or
    181\or 212\or 243\or 273\or 304\or 334\fi
  \advance\countB by -1995
  \multiply\countB by 365
  \advance\countA by \countB
  \countB=\countA
  \divide\countB by 7
  \multiply\countB by 7
  \advance\countA by -\countB
  \advance\countA by 1
  \ifcase\countA\or Sunday\or Monday\or Tuesday\or Wednesday\or
    Thursday\or Friday\or Saturday\fi
\endgroup}

\def\timename{\begingroup
   \countA = \time
   \divide\countA by 60
   \countB = \countA
   \countC = \time
   \multiply\countA by 60
   \advance\countC by -\countA
   \ifnum\countC<10\toks1={0}\else\toks1={}\fi
   \ifnum\countB<12 \toks0={\sevenrm AM}
     \else\toks0={\sevenrm PM}\advance\countB by -12\fi
   \relax\ifnum\countB=0\countB=12\fi
   \hbox{\the\countB:\the\toks1 \the\countC \thinspace \the\toks0}
\endgroup}

\def\timestamp{\dayname, \the\day\ \monthname\ \the\year, \timename}


\def\enma#1{{\ifmmode#1\else$#1$\fi}}

%
\input diagrams.tex
\overfullrule=0pt
\font\tengoth=eufm10  \font\fivegoth=eufm5
\font\sevengoth=eufm7
\newfam\gothfam  \scriptscriptfont\gothfam=\fivegoth 
\textfont\gothfam=\tengoth \scriptfont\gothfam=\sevengoth
\def\goth{\fam\gothfam\tengoth}
%
\font\tenbi=cmmib10  \font\fivebi=cmmib5
\font\sevenbi=cmmib7
\newfam\bifam  \scriptscriptfont\bifam=\fivebi 
\textfont\bifam=\tenbi \scriptfont\bifam=\sevenbi

\font\hd=cmbx10 scaled\magstep1
\def \fix#1 {{\hfill\break \bf (( #1 ))\hfill\break}}
\def \Box {\hfill\hbox{}\nobreak \vrule width 1.6mm height 1.6mm
depth 0mm  \par \goodbreak \smallskip}

\def \coker {\mathop{\rm coker}}

\def \im {\mathop{\rm im}}
\def \Chow {\mathop{\rm Chow}}
\def \length {\mathop{\rm length}}
\def \supp {\mathop{\rm supp}}
\def \deg  {\mathop{\rm deg}}

\def \dimi{{\rm dim}}
\def \dim{\mathop{\rm dim}}
\def \Pic{\mathop{\rm Pic}}

\def \rank {\mathop{\rm rank}}

\def \iso {\cong}
\def \tensor {\otimes}

\def \Hom {{\mathop{\rm Hom}\nolimits}}
\def \hom {{\mathop{\rm Hom}\nolimits}}
\def \GL {{\rm GL}}
\def \Ext {{\rm Ext}}

\def \ext{{\rm Ext}}

\def \h {{\rm h}}

\def \coh{{\rm Coh}}

\def \th {{^{\rm th}}}

\def \B {{\cal B}}
\def \CC {{\bf C}}
\def \C {{\cal C}}
\def \DD  {{\bf D}}
\def \cE  {{\cal E}}
\def \E  {{\cal E}}
\def \F {{\cal F}}
\def \cF {{\cal F}}
\def \FF {{\bf F}}
\def \G {{\cal G}}
\def \GG {{\bf G}}

\def \H {{\rm H}}
\def \cH {{\cal H}}
\def \cI {{\cal I}}

\def \cL {{\cal L}}
\def \L {{\cal L}}
\def \LL {{\bf L}}

\def \O {{\cal O}}
\def \cO {{\cal O}}
\def \P {{\bf P}}
\def \PP {{\bf P}}

\def \RR {{\bf R}}

\def \TT {{\bf T}}
\def \cU {{\cal U}}
\def \U {{\cal U}}
\def \UU {{\bf U}}

\def \ZZ {{\bf Z}}
\def \gm {{\goth m}}

\def\fix#1{\noindent{\bf**** #1 ****}}
\forward {chow complex}{Section}{1}
\forward {Ulrich}{Section}{2}
\forward {atiyah}{Section}{3}
\forward {curves}{Section}{4}
\forward {Projective spaces}{Section}{5}
\forward {surfaces}{Section}{6}
\forward {Appendix}{Section}{7}
\forward{incidence complex}{Theorem}{1.2}
\forward{Chow point}{Theorem}{1.6}
\forward {Segre}{Proposition}{2.6}
\forward {coho of d-uple Ulrich}{Theorem}{4.1}
\forward {chow of Ulrich}{Theorem}{3.1}

\centerline{\hd Resultants and Chow forms via Exterior Syzygies}
\smallskip 
\centerline {\bf David Eisenbud and Frank-Olaf Schreyer
}\footnote{}{\hskip-\parindent
Both authors are grateful for support by the NSF while
accomplishing this research at MSRI. The second author is also 
grateful to the DFG for partial support.  }

\medskip
\centerline{Appendix: {\bf Homomorphisms of some vector bundles on the
Grassmannian}}
\centerline{\bf by Jerzy Weyman}
\smallskip
{\narrower
\noindent{\bf Abstract:} Given a sheaf on projective space $\P^n$
we define a sequence of canonical and easily computable {\it Chow
complexes\/} on the Grassmannians of planes in $\P^n$, generalizing
the well-known Beilinson monad on $\P^n$. If the sheaf has
dimension $k$, then the Chow form of the associated $k$-cycle is the
determinant of the Chow complex on the Grassmannian of planes of
codimension $k+1$.  Using the theory of vector bundles and the
canonical nature of the complexes we are able to give explicit
determinantal and Pfaffian formulas for resultants in some cases where no
polynomial formulas were known. For example, the Horrocks-Mumford
bundle gives rise to a polynomial formula for the resultant of five
homogeneous forms of degree eight in five variables.

}
\medskip
\noindent Let $W$ be a vector space of dimension $n+1$ over a field $K$.
The {\it
Chow divisor\/} of a
$k$-dimensional variety
$X$ in $\P^n=\P(W)$ is the hypersurface, 
in the Grassmannian $\GG = \GG_{k+1}$ of
planes of codimension $k+1$ in $\P^n$, 
consisting of those planes meeting
$X$. The {\it Chow form} is its defining equation.
{}For example the resultant of
$k+1$ forms of degree
$e$ in
$k+1$ variables is the Chow form of $\P^{k}$ embedded by the $e^\th$
Veronese mapping in $\P^n$ with $n={k+e\choose k}-1$. In this paper
we will give a new expression for the Chow divisor,
closely related to Beilinson's monad for sheaves on projective space, 
and
derive new polynomial formulas for Chow forms
in a number of particular cases of the following types:
\smallskip
\noindent {\bf 1. B\'ezout formulas for resultants.}
The classic formula of B\'ezout gives the resultant of
two homogeneous forms in two variables as a determinant of linear forms in
the  Pl\"ucker coordinates of the space generated by the two forms. By analogy we will
call any formula for the Chow form in Pl\"ucker coordinates a {\it B\'ezout
expression} of the Chow form. Our simplest example of a new B\'ezout expression
is for the resultant of three
forms degree 2 in three variables (we also give corresponding formulas for any degree): 
it is the Pfaffian of the alternating matrix of linear forms in the Pl\"ucker coordinates
 $$     
\pmatrix{
 \scriptscriptstyle  0    &\scriptscriptstyle   [245] &
\scriptscriptstyle [345] &\scriptscriptstyle   [135] &
\scriptscriptstyle   [045] &\scriptscriptstyle    [035]
&\scriptscriptstyle    [145] &\scriptscriptstyle     [235]   \cr
\scriptscriptstyle -[245] &  \scriptscriptstyle   0   &\scriptscriptstyle -[235] &\scriptscriptstyle   [035] &  \scriptscriptstyle  [025] &  \scriptscriptstyle  [015] &  \scriptscriptstyle  [125] & \scriptscriptstyle -[125]+[045] \cr
\scriptscriptstyle -[345] &  \scriptscriptstyle [235] &  \scriptscriptstyle  0   & \scriptscriptstyle [134] &\scriptscriptstyle    [035] &   \scriptscriptstyle [034] & \scriptscriptstyle   [135] &     \scriptscriptstyle[234]   \cr
\scriptscriptstyle -[135] & \scriptscriptstyle -[035] & \scriptscriptstyle -[134] &  \scriptscriptstyle  0   &  \scriptscriptstyle [023] &   \scriptscriptstyle  [013] & \scriptscriptstyle [123]-[034]&\scriptscriptstyle   -[123]   \cr
\scriptscriptstyle -[045] & \scriptscriptstyle -[025] & \scriptscriptstyle -[035] & \scriptscriptstyle -[023] & \scriptscriptstyle    0   & \scriptscriptstyle [012] &  \scriptscriptstyle -[015] & \scriptscriptstyle -[024]+[015] \cr
\scriptscriptstyle -[035] & \scriptscriptstyle -[015] & \scriptscriptstyle -[034] &\scriptscriptstyle  -[013] & \scriptscriptstyle  -[012] & \scriptscriptstyle    0  & \scriptscriptstyle [023]-[014] &\scriptscriptstyle  -[023]  \cr
\scriptscriptstyle -[145] & \scriptscriptstyle -[125] & \scriptscriptstyle -[135] &
\scriptscriptstyle-[123]+[034] & \scriptscriptstyle  [015] &
\scriptscriptstyle-[023]+[014] & \scriptscriptstyle 0 &\scriptscriptstyle -[124]+[035] \cr
\scriptscriptstyle -[235] &\scriptscriptstyle [125]-[045] &\scriptscriptstyle -[234] &   \scriptscriptstyle [123] & \scriptscriptstyle [024]-[015] & \scriptscriptstyle [023] & \scriptscriptstyle [124]-[035] &\scriptscriptstyle 0    \cr}
$$\goodbreak
\noindent Here the monomials in the three variables
$a,b,c$ are ordered $a^2,ab,ac,b^2,bc,c^2$
and the brackets
$[ijk]$ denote the corresponding Pl\"ucker coordinates of the net of quadrics. Using 
the theory of rank two vector bundles on $\PP^2$ we can construct many such formulas
for ternary forms of any degree.

\smallskip

\noindent {\bf 2. Stiefel formulas for resultants.}
The Grassmannian is a quotient of an
open set in the variety of $(k+1) \times (n+1)$ matrices over $K$; the
entries of these matrices are called {\it Stiefel coordinates} on the Grassmannian
(or on the Stiefel manifold.) Pulling back the Chow divisor, we get a divisor
whose ideal is generated by a polynomial in the Stiefel coordinates. For example
if $X$ is the rational normal curve this polynomial is the Sylvester determinant.
Even when we cannot express the Chow form of a variety as the determinant or Pfaffian
of a matrix in the Pl\"ucker coordinates, we can sometimes express it as
the determinant or Pfaffian of a map of equivariant vector bundles on the 
Grassmannian. 
Such maps pull back to matrices in the Stiefel coordinates whose determinant
or Pfaffian defines the (closure of the) preimage of the Chow divisor.
We say that such a matrix gives a {\it Stiefel expression\/} for
the Chow form. The classical Sylvester determinant is such
a Stiefel expression. 

Explicit polynomial expressions, in particular Stiefel expressions
for the resultant of $k+1$ forms of degree $d\geq 2$ in $k+1$ variables
(Chow form of the $d$-uple embedding of $\PP^k$) have been known
only for $k\leq 3$ (all $d$) and $k=4, d=2$ (see for example
Gel'fand, Kapranov, and Zelevinsky [1994]. Using our method
and constructions of vector bundles on $\PP^k$
we give new Stiefel expressions. In particular,
the Horrocks-Mumford bundle gives rise to Pfaffian Stiefel expressions
for
the resultants of 5 forms of degrees 4, 6, or 8 in 5 variables.
The matrices involved are too large to exhibit here; but
Macaulay2 programs
for producing them and other new examples can be found at
{\tt http://www.msri.org/****}.

\medskip
We next introduce the basic ideas of this paper, and then describe our main results.
Let
$$
\P^n \lTo^{\pi_1} \FF_l  \rTo^{\pi_2} \GG_l
$$
be the incidence correspondence; that is, let $\FF_l$ be the
set of flags consisting of a point $p \in \P^n$ and an
plane $L\in \GG$ of codimension $l$ in $\P^n$
with $p\in L$. Taking $l=k+1$, the Chow divisor of a reduced irreducible
$k$-dimensional subvariety
$X\subset \P^n$ is  by definition $D_X=\pi_2(\pi_1^{-1}X)$.

If $\F$ is any sheaf whose support is $X$, it follows
that the Chow divisor in $\GG$ is
the codimension 1 part of the support of the sheaf
$\G={\pi_2}_*(\pi_1^{*}\F)$. If $\F$ is generically a vector
bundle of rank $r$ on $X$ then $\G$ will be generically
a vector bundle of rank $r$ on $D_X$. Thus $D_X$ can be recovered
as the codimension 1 part of the (scheme-theoretic) support of $\G$.

If $\G$
happens to be presented by a square matrix with nonzero determinant in the Pl\"ucker
coordinates on $\GG$ (or more generally
by a monomorphism of vector bundles on $\GG$) then the Fitting lemma shows that the
determinant is the $r\th$ power of the 
Chow form of $X$.
One of the central contributions of this paper is to give a simple characterization of
a class of sheaves
$\F$ for which $\G$ has such a presentation: they are the ``Ulrich sheaves''
described below.

More generally, one can use the determinant of a complex,
first introduced (for this purpose!) by Arthur Cayley [1848].
This determinant is in general a rational function, the alternating
product of certain minors in matrices representing the complex.
If $\C$ is a complex of locally free sheaves on $\GG$
whose only homology in codimension 1 is $\H^0(\C)=\G$, then
the determinant of $\C$ is the $r^\th$ power of a form defining
the codimension 1 part of the support of $\G$. Such complexes
were produced from Koszul complexes by Cayley, F.~S.~Macaulay, Jouanolou
and 
other authors who derived expressions for resultants as rational functions
in the Chow or Stiefel coordinates. However, these complexes have
been constructed explicitly for only a limited class of sheaves $\F$. 
For modern results, see
Weyman and Zelevinsky [1994] as well as Jouanolou [1995].
An exposition may be found in the book of
Gel'fand, Kapranov and Zelevinsky [1994]. Of course the Chow form is a polynomial: 
in these rational function expressions the denominator
divides the numerator. However, it is not known how to make the quotient explicit.
Refinements aimed at
reducing the degree of the denominator are an active subject of research;
see for example d'Andrea and Dickenstein [2000].

Grothendieck gave a conceptual framework for these constructions 
in an unpublished letter to David Mumford in 1962;
the details are worked out 
by Knudsen and Mumford in [1976], where the letter is described:
By general theory
there always exist locally free complexes $\C$, well-defined
up to homotopy equivalence, with 
$$
\C \simeq {{\bf R}\pi_2}_*(\pi_1^*\F)
$$ 
in the derived category. If $\F$ is generically a vector
bundle of rank $r$ on $X$, then $\C$ satisfies the conditions
above, and so the determinant of $\C$ is the $r^\th$ power of the
Chow form. We will call such a complex $\C$ a {\it Chow complex for $\F$.}
More
generally, when $\F$ is any coherent sheaf on $\P^n$, whose
support has dimension $k$,
the determinant of a Chow complex for $\F$ is the 
Chow form of the $k$-cycle of $\F$, that is the sum of the
Chow forms of the $k$-dimensional components of the support,
each raised to the power equal to the multiplicity of $\F$
at the generic point of that component.

Our first main result gives a canonical 
Chow complex $\UU_{k+1}(\F)$ for each coherent sheaf $\F$, part of a 
sequence of complexes generalizing the Beilinson monad for $\F$.
The construction is so explicit that it can be made on a computer.
Recall that a plane of codimension $l$ in $\PP^n$ corresponds to
an $(n+1-l)$-quotient of $W$, and thus to an $l$-dimensional subspace of
$W$. We write $U_l$ for the tautological $l$-subbundle on $\GG_l$.

\theorem{first main} For any coherent sheaf $\F$ on $\P^n$
there
is a canonical complex $\UU_l(\F)$ of
vector bundles on $\GG_l$ 
with
$$
\UU_l(\F) \simeq {{\bf R}\pi_{2}}_*(\pi_{1}^*\F).
$$ 
The $e^\th$ term of $\UU_l(\F)$ is
$
\sum_j \H^j(\F(e-j)\otimes\wedge^{j-e}U_l.
$

The complex $\UU_n(\F)$ on $\PP^n$ itself is the Beilinson monad
defined in Eisenbud, Fl\o ystad and Schreyer [2000]. The sheaf $\F$ can be
recovered from $\UU_n(\F)$ simply by taking homology. The sheaf 
$\F$ can be recovered from some
of the other $\UU_l(\F)$ as well: just as one
can recover a variety of dimension
$k$ from its Chow divisor in $\GG_{k+1}$, so one can recover any sheaf $\F$
whose support has dimension at most $k$ from the Chow complex $\UU_l(\F)$
as long as $l>k$. All these matters are explained in \ref{chow complex}

Most significant in our treatment is that we can give
an explicit and canonical description of the maps
in the complex $\UU_l$. Until now, in general,
it has only been possible
to write down the sheaves in such a complex 
(see for example Gel'fand-Kapranov-Zelevinsky [1994] section 3.4E,
``Weyman's complexes"), or to approximate
the maps via a spectral sequence. With
enough vanishing of cohomology it was possible to write
down the maps; but these cases were often not
the ones of primary interest.
Also, previous authors seem only to have considered 
formulas coming from the case where $\F$ is a line bundle.
Our technique allows us to recover  explicit expressions of the Chow
form in all the previously known cases, and, using vector bundles
as in the examples mentioned above, some new ones.

The most useful formulas for the Chow form
occur when the complex $\UU$ has just one nontrivial map $\Psi$
$$
\UU: \cdots\rTo 0\rTo 0\rTo C^{-1}\rTo^\Psi C^{0}\rTo 0\rTo 0\rTo\cdots .
$$
In this case the 
Chow form is the determinant of $\Psi$, and if the bundles 
$C^i$ are direct sums of exterior powers
of the tautological bundles, then one gets a determinantal
expression for the Chow form in Stiefel coordinates.

An even better case occurs when 
$C^{-1}\iso \oplus \O_{\GG}(-1)$, a direct sum of copies of $\O_{\GG}(-1)$,
and 
$C^0$ is a direct sum of copies of $\O_{\GG}$: then the Chow form is given
directly as a determinant in the Pl\"ucker coordinates---that is, we get a 
B\'ezout expression for the Chow form of $\F$ and thus for a
power of the Chow form of the support of $\F$. If $\F$ has
rank 1, or if $\F$ has rank 2 and the map $\Psi$ is skew symmetric
so that we can extract the square root of the determinant as the
Pfaffian, then we get the Chow form of the support of $\F$ itself.

Such cases are considered in \ref{Ulrich}.
Our second main result describes precisely the
conditions on the sheaf $\F$ that are necessary for the Chow complex
$\UU_{k+1}(\F)$ to degenerate to one of these special forms. For example:

\theorem{second main} The Chow complex $\UU_{k+1}(\F)$ above degenerates
to a single map $\O_\GG^d(-1)\to\O_\GG^d$ if and only if the 
module of twisted global sections 
$\oplus_m \H^0(\F(m))$ is a Cohen-Macaulay module
with a linear free resolution.

Here by a linear free resolution we mean a free resolution of the form
$$
\cdots\rTo S^{r_2}(-2)\rTo S^{r_1}(-1)\rTo S^{r_0}.
$$
Such Cohen-Macaulay graded modules
$M$ have been studied by Bernd Ulrich [1984] under the name
``maximally generated maximal Cohen-Macaulay modules" and by
others (see Brennan, Herzog, and Ulrich [1987],
Herzog, Ulrich, and Backelin [1991], and the
references given there) under the names
``linear maximal Cohen-Macaulay modules" or simply ``Ulrich modules"  --- we
shall call the corresponding sheaves {\it Ulrich sheaves\/}. For example,
a line bundle $\F$ on a curve $X$ of genus $g$ embedded in 
$P^n$ is an Ulrich sheaf
if and only if $\F(-1)$ has degree $g-1$ and no global sections;
that is, $\F$ corresponds to a point in $\Pic^{g-1}(X)$ which lies
outside the theta divisor $\Theta \subset \Pic^{g-1}(X)$. 

In \ref{atiyah} we turn to the problem of giving
determinantal and Pfaffian expression for the Chow form of an Ulrich sheaf $\F$.
We can express them
directly in terms of the free resolution of the corresponding module $M$
by using a construction developed in Lejeune-Jalabert and Angeniol [1989]
to describe Atiyah classes. Suppose that
$$
0\to F_c\rTo^{\phi_c}\cdots \rTo F_1\rTo^{\phi_1}F_0
$$
is the linear free resolution of $M$ as above.
Regarding the $\phi_i$ as matrices of elements of $W$, we can compose
them as if they were matrices of linear forms in the exterior
algebra: we write $\Psi_\F:=(1/c!)\phi_1\wedge\phi_2\wedge\cdots\wedge \phi_c$
for this product (defined in a slightly different way in
positive characteristic), which is represented by a matrix of 
forms in $\wedge^cW$. We may identify $\wedge^cW$ with the 
the space of linear forms on $\GG$ and we have:

\theorem{third main} If $\F$ is an Ulrich sheaf, then
$\Psi_\F$ is the (only) nonzero map in the Chow complex $\UU(\F)$.
In particular
the Chow form of $\F$ is $\det \Psi_\F$. If $\F$ is
a vector bundle on a $k$-dimensional variety $X$, and
$\F$ is skew symmetric in an appropriate sense, then
(in characteristic not 2) $\Psi_\F$  is skew-symmetric, 
and the square-root of the Chow form of $\F$ is the
Pfaffian of $\Psi_\F$.

\ref{third main} gives a new method for constructing
resultants and Chow forms:
find Ulrich sheaves (or weakly Ulrich sheaves, or Ulrich
sheaves satisfying the skew-symmetry condition\dots) and then
construct the map $\Psi_\F$.
For this construction one can either use the product formula
above or the definition
of the canonical Chow complex $\UU(\F)$ from maps in
a certain free resolution over the exterior algebra. 

The second part of
this paper gives a number of examples of this method. 
We can be completely explicit in
only a small number of cases, and these sections 
leave open a multitude of theoretical and practical problems.
Central to this pursuit is the

\proclaim Problem. Does every embedded variety $X\subset \PP^n$ have an Ulrich 
sheaf? If $X$ has an Ulrich sheaf, what is the smallest possible rank for such a sheaf?

For example, Brennan, Herzog, and Ulrich [1987] showed that when $X$ is
an arithmetically Cohen-Macaulay curve over an infinite field, or a
a complete intersection, or a linear determinantal variety, then
$X$ has an Ulrich sheaf. Doug Hanes showed in his Thesis under Hochster
[1999] that
the $d$-uple embeddings of $\PP^k$ have Ulrich sheaves when
$k\leq 2$ or $k=3$ and $d=2^r$ is a power of two.

\ref{curves} is devoted to the case of curves. We complete (and reprove) 
the result of
Brennan, Ulrich, and Herzog by showing that every curve in 
$\PP^n$ over an infinite field has skew-symmetric rank 2 Ulrich sheaves.
If the field is algebraically closed it has rank 1 Ulrich sheaves; they are in
one-to-one correspondence with the line bundles of degree $g-1$ having
no sections. Thus there are
B\'ezout expressions for the Chow forms of such curves, generalizing the case
of binary forms ($\PP^1$ and the line bundle $\O_{\PP^1}(-1)$) and the
well-known result that the equation of any plane curve over an algebraically
closed field can be written as the determinant of 
a matrix of linear forms.

Such Ulrich sheaves give rise, in principle, to continuous families of resultant formulas for 
sections of a line bundle on a curve of genus $\geq 1$, but it is not easy to make
such formulas explicit. We illustrate with the case of hyperelliptic
curves, and provide a resultant formula for functions of the form
$a+b\sqrt{f}, \; c+d\sqrt{f}$ where $a,b,c,d$ and $f$ are polynomials
in one variable. We carry out the proof completely only in case the degrees
of the various polynomials are small. In the special case of elliptic curves, we get
a resultant formula for doubly periodic functions written in terms of the
Weierstrass $\wp$-function and its derivative.

We turn in \ref{Projective spaces} to the case where $X\subset \PP^n$ is
the $d^\th$ Veronese ($d$-uple) embedding of $X=\PP^k$. This is the 
case that gives rise to resultant formulas for $k+1$ forms of degree $d$
in $k+1$ variables. We give cohomological criteria for a bundle on
$\PP^k$ to be Ulrich for the $d$-uple embedding.
Following a suggestion of Jerzy Weyman, we use this to show that
every Veronese variety has an Ulrich sheaf, obtained
by applying a certain (unique) Schur functor to the tautological
quotient bundle. This gives a way of writing a power
of the resultant as the determinant of a matrix of linear
forms in the Pl\"ucker coordinates. This may be of computational significance:
the use of resultants in computation is to determine whether or not a set
of polynomials has a common zero; a power of the resultant does this just
as well. 

In this section we also find Ulrich modules of rank 2 for each
Veronese embedding of $\PP^2$. 
We prove a lower bound on the ranks of possible Ulrich modules
and using this and a result of Hartshorne-Hirschowitz on the 
existence of mathematical instanton
bundles we show that rank 2 Ulrich modules exist on the $d$-uple embedding
of $\PP^3$ if and only if $d$ is not divisible by 3. 
On $\PP^4$ we show that the Horrocks-Mumford bundle is weakly Ulrich
for the 4,6, and 8-uple embeddings, and satisfies
the skew-symmetry condition necessary for us to get a Pfaffian Stiefel formula
for the corresponding resultants.

\ref{surfaces} is concerned with the existence of skew-symmetric rank 2 Ulrich sheaves
on various surfaces, and thus with Pfaffian resultant formulas generalizing the B\'ezout
formula for $\PP^2$ given at the beginning of this introduction. We use Mukai's
construction of vector bundles on surfaces, and describe the necessary data. Our main result 
here is the existence of skew-symmetric rank 2 Ulrich modules for certain embeddings
of the plane blown up at a set of points, leading to Pfaffian B\'ezout expressions for the
resultant of 3 ternary forms of degree $d$ with assigned simple base points, valid when
the ideal defining the set of base points is generated in degree  $<d$.

Throughout this paper we rely on a certain construction of
homomorphisms between exterior powers of the tautological bundle
on a Grassmannian, explained in \ref{chow complex}. In \ref{Appendix}
Jerzy Weyman proves---in all characteristics---that 
in fact every homomorphism arises from this construction.

The authors are grateful to 
Hans-Christian v.~Bothmer,
Wolfram Decker, 
Joe Harris,
J\"urgen Herzog,
Michael Kapranov, 
Bernd Sturmfels, 
and Jerzy Weyman 
for discussions of various parts of this  material. Finally,
this paper owes much to experiments made with the computer
algebra system Macaulay2 ; thanks to Dan Grayson and Mike Stillman
for writing it [1993-- ] and for their support in using it for this project.

\goodbreak
\section{chow complex} Chow complexes obtained from the Beilinson monad

As above we
write $\GG_l$ for the Grassmannian of planes of codimension $l$ in $\P:=\PP^n=\P(W)$
and   $\FF_l$ for the flag variety
of flags consisting of a point $p \in \P$ and an
plane $L\in \GG$ of codimension $l$ in $\P$
containing $p$. Throughout this section
we will consider the incidence correspondence
$$
\P \lTo^{\pi_1} \FF_l  \rTo^{\pi_2} \GG_l.
$$

Let $0 \rTo U \rTo W \tensor \O_{\GG_l}  \rTo Q \rTo 0$
be the tautological sequence on the Grassmanian $\GG_l$,
so that $U=U_\l$ is a bundle of rank $l$. We write $E$ for the
exterior algebra $\wedge V$, where $V=W^*$. Any element
$a\in \wedge^p(V)$ induces a homomorphism
$\wedge^pW\to K$ and thus a homomorphism of sheaves
$$
\wedge^pU\hookrightarrow 
\wedge^pW\tensor \O_{\GG_l} \to 
\O_{\GG_l}.
$$
Using the diagonal map 
$\wedge^q U\rTo^{\Delta_U}\wedge^{q-p}U\otimes\wedge^pU$
we get maps 
$$
\wedge^qU\rTo^{(1\tensor a)\Delta_U}\wedge^{q-p}U
$$
for every $p,q$. 

We will use the well-known part a) of the following lemma
heavily. We include part b) for background.

\proposition{Hom computation} Let $U=U_l$ be the tautological
subbundle on $\GG_l$.
\item{a)} The maps above make
$\wedge U$ into a module over $\wedge V$. 
\item{b)}({\bf J. Weyman}) The maps
$$
\wedge^pV \to \Hom(\wedge^q U, \wedge^{q-p}U).
$$
are isomorphisms for all integers 
$p,q$ such that $0\leq q-p,q\leq l$.

\proof a) With notation as above,
the naturality of the diagonal maps shows that the 
diagrams
$$\diagram
\wedge^qU& \rTo& \wedge^qW\otimes \O_{\GG_l}\cr
\dTo^{(1\tensor a)\Delta_U}&& \dTo_{(1\tensor a)\Delta_W}\tensor 1\cr
\wedge^{q-p}U& \rTo&
\wedge^{q-p}W\tensor\O_{\GG_l}
\enddiagram
$$
commute. Since $\wedge W$ is a naturally a module over $E=\wedge V$
by this action
(see for example Eisenbud [1995, Appendix A2.4.1]), so is $\wedge U$.

b) This is proved in an appendix to this
paper by J.~Weyman. In characteristic 0 the result follows
from Bott's vanishing
theorem (see Jantzen [1987]). In arbitrary characteristic it is more delicate.
\Box

We will grade $E$ by the convention that the elements of $V$ have
degree $-1$. As usual we write $E(q)$ for the free graded $E$-module
of rank 1, with generator in degree $-q$. Thus, for example, if $q>p$ then
$\Hom(E(q), E(q-p))=E_{-p}=\wedge^pV$. Recall from 
Eisenbud, Fl\o ystad and Schreyer [2000] that a {\it Tate resolution\/}
is a doubly infinite exact complex of finitely generated free
graded $E$-modules which is minimal in the sense that each free
module maps into $V$ times the next one.
If $\F$ is any
coherent sheaf on $\P$ then there is a 
{\it Tate resolution\/} $\TT(\F)$ naturally
associated to $\F$, which can be computed,
using free resolutions over an exterior algebra, from the 
module of twisted global sections $\oplus_e \H^0\F(e)$. Its
$e^\th$ term is isomorphic to
$$
T^e(\F) = \oplus_j\H^j(\F(e-j))\otimes E(j-e).
$$
For all this see Eisenbud-Fl\o ystad-Schreyer [2000].

We can define an  the additive functor $\UU_l$ from graded free 
modules over $E$ to locally free sheaves on $\GG_{l}$
 by taking
$\UU_l(E(p)) = \wedge^p U$, where $U=U_l$ is the tautological
subbundle, and 
sending a  map $\eta: E(q)\to E(q-p)$ to the map 
$\UU_l(\eta): \wedge^{q} U \to \wedge^{q-p} U$ made from
the element $\wedge^pV$ corresponding to $\eta$.

If $\TT$ is any Tate resolution  then
for $e>>0$ or $e<<0$ we have $\UU_l(T^e)=0$, so
$\UU_l(\F):= \UU_l(\TT)$
is a bounded complex of locally free sheaves
on $\GG_\l.$ 

For example, $\UU_{n}(\F)$ is 
shown by Eisenbud, Fl\o ystad and Schreyer [2001] 
to be a {\it Beilinson
monad\/} for the sheaf $\F$ in the sense that it has the terms above,
and its only homology is $\F$, in degree 0 
(the functor $\UU_{n}$ is called $\Omega$ in that paper).

\theorem{incidence complex}
If $\F$  is a sheaf on $\PP^n$ then the complex
$\UU_l(\F)$ represents $\RR{\pi_2}_*(\pi_1^* \F)$ in the 
derived category of sheaves on the Grassmannian $\GG_l$.

\proof By Theorem 6.1 of Eisenbud, Fl\o ystad and Schreyer [2000],
$\UU_n(\F)$ represents $\F$ in $D^b(\coh(\P^n))$.
We will show first that $\UU_l(\F)={\pi_2}_*(\pi_1^* \UU_n(\F))$,
and second that $\RR^i{\pi_2}_*(\pi_1^*(\wedge^pU_n))=0$ for $i>0$.
It follows 
that 
${\RR\pi_2}_*\pi_1^* \UU_n(\F)\cong {\pi_2}_*(\pi_1^* \UU_n(\F))
=\UU_l(\F)$,
as desired.

On $\FF$ we have inclusions of the universal subbundles
$$
\pi_2^*(U_l)\subset \pi_1^*(U_n)\subset W\otimes \O_\FF.
$$
Pushing the left hand inclusion forward we get a canonical map
$U_l={\pi_2}_*\pi_2^*U_l \to {\pi_2}_*\pi_1^*U_n$, 
and we deduce similar maps on the exterior powers.
To show that these are isomorphisms we may compute fiber by
fiber.
If $u \in \GG_l$ then we will also write
$u \subset W$ for the corresponding $l$-dimensional linear subspace.

Setting $\P'=\P(W/u)\subset \P(W)$,
we have the decomposition
$$
\wedge^pU_n \mid_{\P'} \cong
\oplus_{i=0}^p
\wedge^iu
\otimes
\wedge^{p-i}U'_{n-l},
$$
where $U'_{n-l}$ denotes the tautological sub-bundle on $\P'$.
Thus the map $\wedge^pu\to \H^0(\wedge^pU_n\mid_{\P'})$
is an isomorphism, and all other cohomology of
$\wedge^pU_n\mid_{\P'}$ vanishes.

{}From the base change theorem
(Hartshorne [1977], III,12)  It follows that
$\RR^i{\pi_2}_*(\pi_1^*\wedge^pU_n)=0$ for $i>0$ and
${\pi_2}_*(\pi_1^*\wedge^pU_n) \cong \wedge^p U_l$.
%
%
\Box 

The sheaf $\F$ is determined from $\UU_n(\F)$, the Beilinson
monad, by the formula $\F=\H^0(\UU_n(\F))$. In general we have:

\proposition{determine F} If $\F$ is a coherent sheaf of 
dimension $k$ on $\P$ and $l>k$ then $\F$ is determined
by the complex $\UU_l(\F)$.

\proof
The Tate resolution $\TT(\F)$ is determined by any differential 
$\phi_i: T^i(\F)\to T^{i+1}(\F)$, because $\TT^{\ge i+1}(\F)$ is the minimal 
injective resolutuion of $\im \phi_i$ and $\TT^{\le i}(\F)$ the minimal projective resolution of  $\im \phi_i$. Moreover $\TT(\F)$ determines the Beilinson monad and
hence $\F$. Thus it suffices to reconstruct one of the differentials.

The degrees of the generators of the free module in $T^e(\F)$
range (potentially) from $e-k$ to $e$. Thus the degrees of the
generators of $T^{-1}(\F)$ and $T^{0}(\F)$ range at most from 
$-k-1$ to 0. It follows that these free modules can be
recovered from $\UU(\F)=\UU(\TT(\F)) $. 
By \ref{Hom computation} the map between 
$T^{-1}(\F)$ and $T^0(\F)$  can recovered from $\UU(\F)$ as well.\Box

Now we come to the case involved in the Chow divisor.
Given a finite complex of locally free sheaves on  scheme
$$\B: 0 \rTo \ldots \rTo \B^j \rTo \B^{j+1} \rTo \ldots \rTo 0$$
its determinant bundle is defined as
$$
\det(\B) = \prod_{j \hbox{ even }} \det(\B^j) 
\tensor
\prod_{j \hbox{ odd }} \det(\B^j)^*.
$$ 
If $\B$ is generically exact, then there is a
Cartier divisor called the determinant divisor of $\B$
which measures the part of the homology of $\B$ supported in
codimension 1; 
see Knudsen and Mumford [1976] or Gel'fand, Kapranov and 
Zelevinsky [1994 Appendix A] for the general definition.
If $\F$ is a coherent sheaf on $\P(W)$
with support of dimension $k$, then we define the
Chow divisor of $\F$
to be the usual Chow divisor of the 
$k$-dimensional cycle associated to $\F$---the sum of the 
Chow divisors of the $k$-dimensional components of the support of
$\F$, each with multiplicity equal to the multiplicity of $\F$
on that component. The Chow form
$\Chow(\F)$ is the equation of that divisor; it is a section of
$\O_{\GG_{k+1}}(\deg\F)$ defined up to multiplication by a scalar.
The following  Theorem is a more explicit version the main result of
Knudsen and Mumford [1976] Chapter II. 

\theorem{chow point} Let $\F$ be a coherent sheaf
on $\P(W)$. If $\dim \F = k$ then the Chow divisor of
$\F$ is the determinant divisor of the complex $\UU_{k+1}(\F)$.
Moreover, in codimension 1 the only homology of this
complex is at the $0^\th$ term.

We give a proof for the
reader's convenience: 

\proof We may assume that the ground field is algebraically
closed. 
Since $\UU(\F)$ represents $\RR\pi_2(\pi_1^*\F)$ its divisor
does not pass through any point $u$ of the Grassmannian 
such that $\supp(\F) \cap \P(W/u)= \emptyset$. For a general point $u$
of a component of the zero locus of  $\Chow(\F)$
the subspace $\P(W/u)$ meets the support of $\F$ in a single
point which belongs to a unique component $X$ of the support.
We have
$$
\dimi_{\kappa(u)} ({\pi_2}_*\pi_1^* \F) \tensor \kappa(u) 
= \dimi_{\kappa(u)} \H^0(\F \tensor \O_{\P(W/U)})
= \length(\F \tensor \O_{\P(W/u),X})
$$
and the higher direct images vanish.
\Box

\section{Ulrich} Ulrich Sheaves

If $\UU(\cF)$ is a two term complex then the determinant 
section of $\UU(\cF)$ is is the
determinant of a morphism between bundles.
This situation corresponds to the case where the Tate resolution
of $\F$ has ``betti diagram" of the form:
$$
\matrix{ 
h^k\cF(-k-3)  &h^k\cF(-k-2) & h^k\cF(-k-1) & h^k\cF(-k) & 0&0& \cr
       0 & 0 &  h^{k-1}\cF(-k) & h^{k-1}\cF(-k+1) &0&0 \cr
                  \vdots&\vdots&\vdots&\vdots &\vdots&\vdots \cr
0&0&h^1\cF(-2)&h^1\cF(-1)&0&0\cr
0&0&h^0\cF(-1)&h^0\cF&h^0\cF(1) &h^0\cF(2)&  \cr}
$$
Here, by the betti diagram of $\TT(\F)$ we mean
the table whose $(i,j)$ entry
is the number of generators
of degree $j-i$ required by the $j^\th$ free module
$T^i$ in $\TT(\F)$---by 
Eisenbud-Fl\o ystad-Schreyer [2000] this is the dimension of $\H^i(\F(j-i))$.
(This is almost the same as the betti diagram in the programs
Macaulay of Bayer and Stillman, or Macaulay2 of Grayson and Stillman,
except that we think of the arrows in the resolution as going
from left to right. This change of convention is convenient because
of the fact that the generators of $E$ have negative degree.)

For reasons that will become clear in a moment, we will call
a sheaf $\F$ with cohomology as above a {\it weakly Ulrich sheaf\/}.

An even better situation occurs when the Tate resolution has 
betti diagram of the form
$$
\matrix{ 
\ldots&h^k\cF(-k-3)  &h^k\cF(-k-2) & h^k\cF(-k-1) &0 & 0&0& \cr
      & 0 & 0 &  0 & 0 &0&0 \cr
                  &\vdots&\vdots&\vdots&\vdots &\vdots&\vdots \cr
&0&0&0&0&0&0\cr
&0&0&0&h^0\cF&h^0\cF(1) &h^0\cF(2)&\ldots }
$$
\goodbreak 
\noindent In this case we see from the previous section that 
the Chow form of $\cF$ is the determinant of the $h^0(\F)
\times h^k\F(-k-1)$ matrix whose entries are
linear forms in the 
Pl\"ucker coordinates on the Grassmannian $\GG_{k+1}$.
(It follows that $h^0\cF = h^k\cF(-k-1)=\deg(\F)$,
which one can easily see in other ways as well.)

Modules whose associated sheaf have this sort of
Tate resolution were first studied by Bernd Ulrich in [1984]. 
We will call them {\it Ulrich sheaves.\/} Thus a $k$-dimensional sheaf $\F$ 
on $\P$ is
an Ulrich sheaf  if $\F$ has no intermediate cohomology---that 
is, $\H^q(\F(d))=0$ for
$1\leq q\leq k-1$ and all $d$---and
$\H^0(\F(j))= 0$ for $j<0$ while
$\H^k(\F(j))= 0$ for  $j\geq -k$.
Since an Ulrich sheaf has no intermediate cohomology,
its restriction to the nonsingular part of $X$ is automatically
a vector bundle.

We can characterize Ulrich sheaves without referring to all the cohomology
in several elementary ways. Since every 0-dimensional sheaf
is an Ulrich sheaf, we will henceforward ignore this case.

\proposition{Ulrich char} Let $\F$ be a coherent, $k$-dimensional
sheaf on the projective space $\P=\P^n$ over $K$ with $k>0$. The following
are equivalent:
\item{a)} $\F$ is an Ulrich sheaf.
\item{b)} $\H^i \F(-i)=0$ for $i>0$ and $\H^i \F(-i-1)=0$ for $i<k$. 
\item{c)} If the support of $\F$ is a scheme $X$, then for 
some (respectively all) finite linear projections $\pi:X\to \P^k$
the sheaf $\pi_*\F$ is the trivial sheaf $\O_{\P^k}^t$ for some $t$.  
\item{d)} The
module $M:=\H^0_*(\F):=\oplus_d \H^0(\F(d))$ of twisted global sections 
is an {\it Ulrich module,\/} in the sense of 
Backelin and Herzog [1987]; that is, $M$ is a
Cohen-Macaulay module of dimension $k+1$ over the homogeneous coordinate ring
$S=k[x_0,\cdots,x_n]$ of\/ $\P$, whose number of generators is
equal to $\deg \F$, or equivalently whose $S$-free
resolution
$$ 
\FF: 0 \rTo F_{n-k} 
\rTo^{\varphi_{n-k}} \ldots 
\rTo^{\varphi_2}  F_1 
\rTo^{\varphi_1}  F_0 
\rTo M \rTo 0
$$
is linear
in the sense that $F_i$ is generated in degree $i$ for every $i$.

\proof a) $\Rightarrow$ b) is trivial.

b) $\Rightarrow$ c): By the finiteness and linearity of $\pi$
we have $\H^i(\F(j))=\H^i((\pi_*)\F(j)$.
The vanishing of cohomology of b) gives vanishing for $\pi_* \F$ which
characterizes the trivial vector bundles on $\P^k$.
 
c) $\Rightarrow$ d): By $c)$  $M=\H_*^0(\F)$ is a free
module over $K[x_0,\ldots,x_k]=\H_*^0(\O_{\PP^k})$ generated in degree 0.
Thus $M$ is a linear Cohen-Macaulay module, that is an Ulrich
module.

d) $\Rightarrow$ a): The equivalence of the two
characterizations of Ulrich modules given in d) may be found in
Brennan, Herzog, and Ulrich [1987, Prop.~1.5]. 
A graded $S$-module
$M$ is 0-regular if and only if the free resolution of $M_{\geq 0}$ is linear
is proved in Eisenbud-Goto [1984] (see also Eisenbud [1995, Theorem 20.18]).
If $M$ is a $k+1$-dimensional Cohen-Macaulay module with
linear resolution, then the associated sheaf $\F$ is also 0-regular.
The Cohen-Macaulay property of $M$ gives the vanishing of the
intermediate cohomology of $\F$, and (since $\dim M = k+1 >1$)
also shows that
$M=\H^0_*(\F)$. Thus $\H^0(\F(j))=0$ for $j<0$, and $\F$ is Ulrich.
\Box

{}From the linearity of the resolution $\FF$ of an Ulrich module 
$M$ it follows, for
example, that the rank of
$F_{i}$ is ${n-k\choose i}\cdot\rank F_0$; to see this
reduce modulo a maximal
regular sequence, and observe that $M$ must reduce to a direct
sum of copies of the residue field $K$. In particular,
$\rank F_{n-k}=\rank F_0$, and this rank is
equal to the degree of $\F$. (For more details, see for example
Brennan, Herzog, and Ulrich [1987].) The same kind of argument gives:

\corollary{coho of Ulrich} If $\F$ is an Ulrich sheaf of dimension
$k$ on $\P^r$ then 
$
\chi(\F(e))=h^0(\F){e+k\choose k}.
$\Box

In \ref{coho of d-uple Ulrich} we will generalize this to
sheaves on $X$ that are Ulrich sheaves for the $d$-uple embedding
of $\F$.

In our applications we will be particularly interested in the
case where the Ulrich sheaf is a vector bundle on its support,
and is
self-dual up to a twist. In this case the criterion above
can be simplified:

\corollary{self-dual Ulrich}
Let $\F$ be a vector bundle on a $k$-dimensional Gorenstein scheme
$X\subset \P^r$. If $\F\iso \F^*(k+1)\otimes \omega_X$, then
$\F$ is an Ulrich sheaf on $\P^r$ if and only if $\F$ is 0-regular.

\proof The 0-regularity implies that $\H^i(\F(j))=0$ for
$j>-i$. The rest of the necessary vanishing follows from Serre duality.\Box

Brennan, Herzog and Ulrich discovered in [1987] that 
linear determinantal varieties have rank one Ulrich modules,
so we can
give B\'ezout expressions for their Chow forms using the ideas
above. This series of examples includes rational normal scrolls,
Bordiga-White surfaces and many more.
We can give a different description of their
Ulrich modules as follows:

\example{linear determinantal varieties} Let $\varphi \colon F \to G$
with $F=\oplus_1^f \cO$ and $G=\oplus_1^g \cO(1)$, $f\le g$, a linear 
$f \times g$ matrix on $\PP^n$ which drops rank in expected
codimension $(f-g+1)$. The Eagon-Northcott type complex 
$$
0 \to \Lambda^f F \tensor D_{f-g+1} G^* \to \ldots \to \Lambda^{g} F
\tensor G^* \to \Lambda^{g-1} F \to \cF \to 0,
$$
see Eisenbud [1995, Theorem A2.10], is a linear resolution of a module
annihilated by the maximal minors of $\varphi$ and has
length $f-g+1$. It is thus the resolution
of an Ulrich sheaf
on $X = V(I_g(\varphi))$, and one can check that the sheaf has rank
1 (it is isomorphic, in the generic case, to $I_{g-1}\varphi'$,
the ideal generated by the $g-1\times g-1$ minors of the submatrix
$\varphi'$ obtained from $\varphi$ by omitting one row.) Hence the
Chow form of
$Chow(X)=Chow(\cF)$ is polynomial of degree $f \choose {g-1}$ in the Pl\"ucker
coordinates, and $\deg X = {f \choose {g-1}}$.

\example{cubic scroll} Consider the scroll $S(2,1) \subset \PP^4$ defined by
$$ \varphi=\pmatrix{x_0&x_1&x_3 \cr x_1 & x_2 & x_4 \cr}.$$
Using the Ulrich sheaf $\cF$ as above and 
\ref{chow of Ulrich} we obtain its Chow form as the determinant
of the matrix

$$\pmatrix{ [034] & [013]       & [023] \cr
           -[134] & [023]+[014] & -[123]-[024] \cr
            [234] & -[024]   & [124] \cr}. $$

The Chow forms of  rational normal scrolls have further interpretations:
Consider $(r+1)$-dimensional spaces $\alpha$ of sections of bundles
$\oplus_{i=1}^r \O_{\P^1}(d_i)$. The Chow form of the scroll $S(d_1,\ldots,d_r)
\subset \P^N$ with $N+1=\sum_i (d_i+1)$ describes those $\alpha$, 
where the minors of the 
corresponding morphism
$$ 
\O_{P^1}^{r+1} \rTo^\alpha \oplus_{i=1}^r \O_{\P^1}(d_i)
$$
have a common zero. (Such formulas were also worked out
by Henri Lombardi and J.-P.~Jouanolou 
(unpublished).)
 
In case of
$S(2,1)$ there is also the interpretation for plane conics with one assigned 
base point: Since $S(2,1)$ is the image of $\PP^2$ by the linear system of conics with
a single assigned base point, say $(1:0:0)$, its Chow form describes
those 3-dimensional subspaces of conics which have a further base point.
 In section 5 we will generalize this examples to forms of any degree on 
$\PP^2$ with several simple assigned base points.

{}From the point of view of examples, it is interesting to note that if
two schemes in projective spaces support (weakly) Ulrich sheaves, then
so does their Segre product:

\proposition{Segre} Let $\F_1$ be a coherent sheaf on $\P(W_1)$ and let
$\F_2$ be a coherent sheaf on $\P(W_2)$. Set
$d= \dim(\F_1).$ Let $\G$ be the
Segre product of $\F_1$ with $F_2(d)$ on
$\P=\P(W_1\otimes W_2)$; that is, 
$\G=(\pi_1^*\F_1)\otimes(\pi_2^*\F_2(d))$ on the Segre variety
$\P(W_1)\times \P(W_2)\subset \P$.
\item{a)} If $\F_1, \F_2$ are weakly Ulrich, then $\G$ is weakly Ulrich.
\item{b)} If $\F_1, \F_2$ are Ulrich, then $\G$ is Ulrich.

Of course a similar result holds for the Segre product
of $\F_1(\dim \F_2)$ and $\F_2$.

\proof Both parts follow easily from the K\"unneth formula
$$
\H^i(\G(d))=\oplus_{i=j+k}H^j(\F_1(d))\otimes \H^k(\F_2(d)).
$$
For example, in part a) we need 
$\H^j(\F_1(-j-k-2))\otimes \H^k(\F_2(d-j-k-2))=0$
when $j+k<d+\dim F_2$. If $j<d$ then the first factor
vanishes since $\F_1$ is weakly Ulrich, while if $j=d$
then the second factor vanishes for the same reason.\Box

\corollary{product formula} With notation as in \ref{Segre}, 
suppose that 
$\F_1, \F_2$ are Ulrich of dimensions $d_1, d_2$.
The map over the exterior
algebra 
$E_{\rm Segre}=\wedge((W_1\otimes W_2)^*)=
\wedge(W_1^*\otimes W_2^*)$
$$
\H^{d_1+d_2}(\G(-d_1-d_2-1))\otimes\omega_{E_{\rm Segre}}
\rTo\H^0(\G)\otimes \omega_{E_{\rm Segre}}
$$
is derived from the tensor product of the
corresponding maps for $\F_1$ and $\F_2$
over $\wedge W_1^*$ and $\wedge W_2^*$ respectively via the
canonical injection 
$\wedge W_1^* \otimes \wedge W_2^*\subset \wedge(W_1^*\otimes W_2^*)$.
\Box

It follows that in situations where
we can compute a B\'ezout expression for the Chow forms of
$\F_1$ and $\F_2$ we can also
compute a B\'ezout expression for the Chow form of the Segre
product. Similar remarks and formulas hold in the case of 
weakly Ulrich sheaves and Stiefel expressions of the Chow form.

\section{atiyah} Chow forms as determinants and Pfaffians

Throughout this section we will work with a sheaf $\F$ of
dimension $k$ on $\PP^n=\P(W)$. For simplicity, we will write
$\UU$ for the functor $\UU_{k+1}$ defined in \ref{chow complex}.
We set $c=n-k$, the codimension of $\F$.

As we have seen in the previous section,
if $\F$ is weakly Ulrich
then the complex $\UU(\F)$ is given by a single map
$\Psi_\F: \UU(T^{-1}\F)\to \UU(T^0)\F$
of vector bundles on the Grassmannian, and the Chow form
of $\F$ is the determinant of $\Psi_F$. In this section we
will make $\Psi_F$ explicit. The tools we develop
will allow us to show that if $\F$ is skew symmetrically self-dual
in a natural sense then the complex $\UU(\F)$ is skew
symmetric, and in particular $\Psi_\F$ is skew symmetric.
When $\F$ is also weakly Ulrich, the square root
of the the Chow form of $\F$ is the Pfaffian of $\Psi_\F$.
In particular, when $\F$ is in addition a sheaf of rank
2 supported on a variety of $X$, the Chow form of $X$
itself is the Pfaffian of $\Psi_\F$.

The matter is simplest in the Ulrich case, and
we describe this first:
Suppose $\F=\widetilde M$ is an Ulrich sheaf
and let $c=n-k$
be its codimension. By \ref{Ulrich char} $M$
has a linear free resolution
$$ 
L:  
0\rTo  L^{-c} \rTo^{\alpha_{-c}}
\cdots 
\rTo  L^{-1}
\rTo^{\alpha_{-1}}  L^0,
$$
Where $L^{-i}=S\otimes P_i$ for some 
finite dimensional vector space $P_i$ 
concentrated in degree $i$. Each map $\alpha_i$
in the resolution 
corresponds to a map $P_{i}\to W\otimes P_{i-1}$,
and because the maps of free $S$-modules compose to 0,
the composite 
$P_{i}\to W\otimes P_{i-1}\to W\otimes W\otimes P_{i-2}$
has image contained in $\wedge^2W\otimes P_{i-2}$. More generally,
each composite
$
P_{i}\to (\otimes_1^jW)\otimes P_{i-j}
$
has image in $\wedge^jW\otimes P_{i-j}$.

In particular we get a map 
$\varphi_{-c,0}: P_c\to \wedge^cW\otimes P_0$.
We may identify
$\wedge^cW$ with the space of linear forms on the
Grassmannian $\GG=\GG_{k+1}$ of planes of codimension
$k+1$, and thus $\varphi_{-c,0}$ gives a map
$\Psi_\F: \O_\GG(-1)\otimes P_c \to \O_\GG\otimes P_0$. 

\theorem{chow of Ulrich} If  $\F=\widetilde M$ is
an Ulrich sheaf on $\P(W)$, then with notation as
above, the ranks of $P_c$ and $P_0$ are the same
and the Chow form of $\F$ is the determinant of the 
map $\Psi_\F$.\Box

\noindent{\sl Remark:\/} 
Angeniol and Lejeune-Jalabert [1989] define maps
generalizing the
$
P_{i+1}\to \wedge^jW\otimes P_{i-j}\subset(\otimes_1^jW)\otimes P_{i-j}
$
for the free resolution of any module $M$ and use them
to construct the Atiyah classes of $M$.
Because each $L^i$ is a free module generated in a single degree,
our situation is
simpler. In particular, in our case 
the maps themselves---not just their
cohomology classes---are well-defined.
\smallskip

To generalize \ref{chow of Ulrich} we replace 
free resolutions by the
linear free monads studied in 
Eisenbud, Fl\o ystad and Schreyer [2001].
Here is a summary of the theory (references
to EFS refer to that paper):
Let $\F$ be any coherent sheaf on $\PP(W)$.
By by EFS, Example 8.5 and Proposition 8.6, there
is a unique complex
$$
L=\LL(\F):\quad \cdots 
\rTo^{\alpha_{-2}} L^{-1}\rTo^{\alpha_{-1}}
L^0\rTo^{\alpha_0} L^1\rTo^{\alpha_1}
\cdots 
$$
such that  $L^{-i}=S(-i)\otimes P_i$
and $L^i=0$ if $|i|>n$,
with the property that the sheafification of the 
homology of $L$ is zero except for $\widetilde{\H^0(L)} = \F$.
The complex $\LL(\F)$ is called the linear free
monad of $\F$. It is functorial in $\F$, and may be
constructed from the Tate resolution 
$$
\TT(\F): \cdots\to T^{-1}(\F)\to T^0(\F)\to T^1\to\cdots;
$$
in fact, 
$\LL(\F)=\LL(P)$, the complex corresponding to 
the graded $E$-modules $P:=\im (T^{-1}(\F)\to T^0(\F))$ under
the Bernstein-Gel'fand-Gel'fand correspondence. Thus
$\LL(F)$ has the form
$$
\cdots 
\rTo P_1\otimes S= L^{-1}
\rTo P_0\otimes S= L^0
\rTo P_{-1}\otimes S= L^1
\rTo 
\cdots.
$$
For example, the linear
free monad for an Ulrich sheaf is equal to the 
minimal free resolution. 

Associated to $L=\LL(\F)$ are the maps
$$
\varphi_{-i,-j}:\quad P_i\to\wedge^{i-j}W\otimes P_j.
$$
adjoint to the multiplication map $\wedge^{i-j}V\otimes P_i \to P_{i-j}$
that defines the $E$-module structure on $P$. These may also
be computed from the differentials of $L$, as above.

We can now describe the ``middle'' map $T^{-1}(\F)\to T^0(\F)$ in
the Tate resolution $\TT(\F)$.
Let $\gm$ be the ideal of elements of negative
degree in $E$ (the augmentation ideal) and define graded
vector spaces $A$ and $B$ by
$$
A=P/\gm P\qquad B=\{p\in P\mid \gm p = 0\}. 
$$
A projective
cover $F\to P$ is a minimal map from a free $E$-module $F$
onto $P$. It follows from Nakayama's Lemma that $F\iso E\otimes A$.
A projective
cover is determined by the data of a splitting $\eta: A\to P$ (as 
graded vector spaces) of the natural projection map $P\to A$.
Dually, an injective envelope $P\to G$ is uniquely determined
by a splitting $\pi: P\to B$ of the inclusion $B\subset P$;
we take $G\iso \omega_E\otimes B$, and the map from $P$
is the unique map to $\omega_E\otimes B$ whose composition
with the projection to $(\omega_E)_0\otimes B = B$ is $\pi$.

The composition of $F\to P$ and $P\to G$ is a map
$$
\varphi_P: E\otimes A\to \omega_E\otimes B
$$
whose image is $P$. 
We define $\Psi_\F=\UU(\varphi_P)$, which
is a map of vector bundles
on the Grassmanian $\GG_{k+1}$.
For example if $L$ is a free resolution of
an Ulrich sheaf then, by Eisenbud-Fl\o ystad-Schreyer [2001], Proposition 8.7,
$A=P_c$ and $B=P_0$, so in that case
$\varphi_P$ is the map induced by the map $\varphi_{-c,0}$
defined before \ref{chow of Ulrich}, and the map $\Psi_\F$
is the same as the one given there.
(No choice of $\eta$ and $\pi$ is involved because 
$A= P_{-c},\; B=P_0$ in that case.)

\theorem{chow form of weakly Ulrich}
If $\F$ is a weakly Ulrich sheaf of dimension $k$ on $\P(W)$,
with linear monad $\LL(P)$, 
then the Chow form of $\F$ is $\det \UU(\Psi_\F)$.

\proof  By \ref{chow point} the Chow form of $\F$ is
the determinant of the complex
$\UU\TT(\F)$. Since $\F$ is weakly Ulrich, this complex
consists of a single map:
$$
\UU(\F)=\UU(T^{-1}\to T^0) = \UU(\varphi_P)=\Psi_\F.
$$
\Box

\goodbreak
\bigskip
\noindent{\bf The skew symmetry of $\UU(\F)$}
\medskip

We now show that appropriate symmetry or skew symmetry of $\F$
makes $\UU(\F)$ symmetric or skew symmetric.
The functor
$$
D: \F \mapsto \E xt^c(\F,\omega_{\PP^n})(k+1)
$$ 
defines a duality on the category of $k$-dimensional Cohen-Macaulay sheaves
on $\PP^n$ and there is 
a canonical morphism $\iota: \F\to \DD\DD(\F)$.
Let $\epsilon=\pm 1$. 
As with any duality, we
say that a morphism   
$\sigma: \F \rTo D(\F)$ is
$\epsilon$-symmetric if 
$$\diagram
&&DD(\F)\cr
&\ruTo^\iota\cr
\F&\epsilon&\dTo_{D\sigma}\cr
&\rdTo_\sigma\cr
&&D(\F)
\enddiagram
$$
commutes up to the sign $\epsilon.$ 
In case $\epsilon=1$ we say that $\F$ is symmetric;
if $\epsilon=-1$ then $\F$ is called skew symmetric.

\theorem{symmetry}
Suppose that $\F$ is a
Cohen-Macaulay sheaf of dimension $k$
on $\P^n$. Any
$\epsilon$-symmetric isomorphism $\F\to D(\F)$,
induces an $\epsilon$-symmetric isomorphism
$$
\UU(\F)\to \Hom_{\GG_{k+1}}(\UU(\F), \O_{\GG_{k+1}}(-1))[1].
$$ 
In particular the map
$\UU T^{-1}(\F) \rTo^{\Psi_\F} \UU T^0(\F)$
is $\epsilon$-symmetric, and for $j>1$ the map
$\UU T^{-j}(\F) \rTo \UU T^{-j+1}(\F)$
is dual to 
$\UU T^{j-1}(\F) \rTo \UU T^{j}(\F)$.

If $\F$ is skew symmetric we define the 
Pfaffian of the skew symmetric complex $\UU(\F)$ 
by taking an appropriate Pfaffian of the middle map 
$\UU T^{-1}(\F) \rTo^{\Psi_\F} \UU T^{0}(\F)$
times the alternating product of those terms from the 
determinant of $\UU(\F)$ that are associated with the maps
$\UU T^{-j}(\F) \rTo \UU T^{-j+1}(\F)$
for $j>0$. The determinant of $\UU(\F)$ is then 
the square of the Pfaffian of $\UU(\F)$. 

\corollary{det as pfaff} Assume that the
characteristic of the ground field is not 2.
If $\F$ is a skew-symmetric
Cohen-Macaulay sheaf of rank $2$
on a $k$-dimensional subscheme $X\subset \PP^n$
such that $\wedge^2\F \iso \omega_X(k+1)$, then the
Chow form of $X$ is the Pfaffian of the complex $\UU(\F)$.
In particular
if $\F$ is weakly Ulrich, then
the Chow form of $X$ is the Pfaffian of the skew-symmetric
map of vector bundles $\UU(\Psi_\F).$

\noindent{\sl Remark.\/} In order to include the case
of characteristic 2 we would have to add the condition
that the duality $D$ is alternating, not just skew symmetric,
and then prove the corresponding result for $\Psi_\F$.
We leave this task to the interested reader.
\smallskip

\noindent{\sl Proof of \ref{det as pfaff}.\/}
The skew-symmetric pairing 
$\F\otimes \F\to \wedge^2\F\iso \omega_X(k+1)$
gives rise to a skew symmetric isomorphism 
$$
\F\to \Hom(\F,\omega_X(k+1)\iso 
\Ext^c(\F,\omega_{\PP^n})\iso \DD(\F).
$$
The rest follows from \ref{symmetry} and the discussion above.
\qed

To prove \ref{symmetry} we will first analyze the map
on linear monads induced by the (skew) symmetric
isomorphism $\F\to D\F$. From this analysis will
come a certain symmetry property of $\varphi_P$.
The map $\varphi_P$ may be represented by a matrix
of elements of $\wedge W$.
An element $\alpha\in \wedge^{t} W$ induces
for any integer $s$ a map
(which we again call $\alpha$) defined by
$$
\wedge^{s}U\otimes \wedge^vV \rTo^\alpha \wedge^{s+t-v}U\qquad:\qquad
u\otimes e \mapsto \alpha(e)(u)
$$
where $\alpha(e)\in \wedge^{v-t}V$ acts on $U$ as described
in \ref{chow complex}. The map $\Psi_\F$ is constructed
from these pieces, so we will derive a symmetry property
for $\Psi_F$.

The difficulty of proving \ref{symmetry} comes
from the delicacy of the signs involved. For example,
consider the case where the map 
$\varphi_P: E(k+1-i)\to \omega_E(i)$ in 
\ref{chow form of weakly Ulrich} is given by a 
$1\times 1$ matrix whose entry is in $\wedge^{c+2i}W$.
One might suppose that any $1\times 1$ matrix would
be symmetric, and correspond 
to a symmetric map of vector bundles
$\Psi: (\wedge^iU)^*\iso\wedge^{k+1-i}U\otimes \wedge^vV\to \wedge^iU$ on the 
Grassmannian. But actually $\Psi$ is symmetric 
if $i$ is even and skew symmetric if $i$ is odd.
The general result we need is the following:

\lemma{gr signs}
Set $c=v-k-1$ and let $\alpha\in \wedge^{c+i+j}W$.
The dual into $\O_{\GG_{k+1}}(-1)$ of the map
$$
\wedge^{k+1-i}U\otimes \wedge^vV \rTo^\alpha \wedge^jU
$$
is the map
$$
\wedge^{k+1-j}U\otimes \wedge^vV 
\rTo^{(-1)^{k(i+j)+ij}\ \cdot\ \alpha}
\wedge^iU.
$$

\noindent{\sl Proof of \ref{gr signs}.\/}
We identify 
$\wedge^iU$ with $\Hom(\wedge^{k+1-i}U\otimes \wedge^vV,\O_{\GG_{k+1}}(-1))$
via the map $\tau$ sending $\beta\otimes e\in \wedge^{k+1-i}U\otimes \wedge^vV$
to the functional
$$
\tau:\ \wedge^iU \ni\chi\mapsto (\chi\wedge\beta)(e)\in \O_{\GG_{k+1}}(-1).
$$
We must show that the diagram
$$
\diagram
\wedge^{k+1-i}U\otimes \wedge^vV & \rTo^\alpha & \wedge^j U\cr
\dTo^\tau&&\dTo_{\tau^*}\cr
(\wedge^i U)^* &  \rTo_{\alpha^*} &(\wedge^{k+1-j}U\otimes \wedge^vV)^*
\enddiagram
$$
commutes up to a sign of $(-1)^{k(i+j)+ij}$. Although this
is a diagram of vector bundles, we may reduce the
problem to one of vector spaces by working fiberwise. 
For each $p\in \GG_{k+1}$ the
fiber $U_p$ of $U$ is a subspace of $W$, and
the action of $V$ on $U_p$ is induced by
its action on $W$. Thus the annihilator of $U_p$ in $V$ acts
as zero on $U_p$, and we may therefore replace $W$ by $U_p$
and $V$ by $U_p^*$, and assume that  $U=W$
so that $v=k+1$ and $c=0$.

{}From the definitions we see that 
$$\eqalign{
\alpha^*\tau(\beta\otimes e):\quad &\gamma\otimes e \mapsto 
\bigl[\bigl((\alpha(e))(\gamma)\bigr)\wedge \beta\bigr](e), \cr
\tau^*\alpha(\beta\otimes e):\quad &\gamma\otimes e \mapsto 
\bigl[\bigl((\alpha(e))(\beta)\bigr)\wedge\gamma\bigr](e).
}$$
Since these expressions are multilinear in $\alpha, \beta, \gamma$
it suffices to check the case where $\alpha, \beta, \gamma$ are
products of elements in some fixed basis 
$\{x_1,\dots,x_v\}$ of $W$.
Set
$a=\alpha(e)\in\wedge^{k+1-i-j}V$.
The expressions are both zero unless
$a(\beta)\wedge \gamma$ is a scalar times the product of 
all the basis elements $\{x_1,\dots,x_v\}$. Under this 
assumption, what we are trying to prove is equivalent
to the statement that
$
a(\gamma)\wedge \beta = (-1)^{k(i+j)+ij} a(\beta)\wedge \gamma.
$

Let $\overline \alpha$ be the
element of $\wedge^{k+1-i-j}W$ such that
$\alpha(e)(\overline \alpha) = a(\overline \alpha)=1$. 
Our assumptions imply that we can
can factorize $\gamma$ and $\beta$ 
as
$\gamma=\gamma'{\overline\alpha}$ and
$\beta=\overline\alpha \beta'$.
With this notation
$$\eqalign{
a(\gamma)\wedge\beta&=
a(\gamma'\wedge\overline \alpha)\wedge \beta =
(-1)^{\gamma'a}\gamma'\wedge \beta =
(-1)^{\gamma'a}\gamma'\wedge \overline \alpha \wedge \beta',\cr
a(\beta)\wedge\gamma&=
a(\overline \alpha\wedge\beta')\wedge \gamma =
\beta'\wedge \gamma =
\beta'\wedge \gamma'\wedge \overline \alpha  =
(-1)^{(\gamma'+a)\beta'} \gamma'\wedge \overline \alpha \wedge\beta',
}$$
where we have also written $\gamma', a$, and $\beta'$ for the
degrees of these elements.
Thus the diagram commutes up to the sign
$
(-1)^{(\gamma'+a)\beta'+\gamma'a}.
$
But
$
(\gamma'+a)\beta'+\gamma'a=
\gamma(\beta-a)+(\gamma-a)a=
\gamma\beta-a^2=(k+1-i)(k+1-j)+(k+1-i-j)^2
$
and this is congruent modulo 2 to
$k(i+j)+ij$
as required.
\qed
\medskip

\noindent{\sl Proof of \ref{symmetry}}.
We will show
that the
``middle'' differential
$$\UU T^{-1}(\F) \rTo^{\Psi_\F} \UU T^0(\F)$$
is $\epsilon$ symmetric. This condition depends
on an identification of 
$\UU T^0(\F)$ with the dual of $\UU T^{-1}(\F)$. 
Changing this identification is the same as 
multiplying $\Psi_\F$ by an automorphism of
its source or target, so it suffices to show
that $\Psi_\F$ times such an isomorphism is
$\epsilon$ symmetric.

Once we know that the middle differential is
$\epsilon$ symmetric,
we can take the injective resolution
of $P$, from which the positively indexed maps of
$\UU(\F)$ are made, to be dual to the free resolution
of $\P$ from which the negatively indexed maps of 
$\UU(\F)$ are made. 

To analyze $\Psi_\F$ we will make use of the analysis
of $\varphi_P$ described before \ref{chow form of weakly Ulrich}.
We have decompositions
$$
\eqalign{
T^{-1}(\F)&=\sum_i A_{c+i}\otimes E(-c-i)\ {\rm and}\cr
T^{0}(\F)&=\sum_j B_{-j}\otimes \omega_E(j).
}$$
In terms of this decomposition, the $(i,j)$ component
of $\varphi_P$ is $\pi_j\varphi_{-c-i,j}\eta_{-c-i}$. 
Applying the functor $\UU$
we see that $\Psi_\F$ decomposes into maps
$$\eqalign{
\UU(A_{c+i}\otimes E(-c-i)) = 
\UU(A_{c+i}\otimes \wedge^vV \otimes \omega_E(k+1-i)) = 
&A_{c+i}\otimes \wedge^vV \otimes \wedge^{k+1-i}U\cr
\rTo^{(\Psi_\F)_{i,j}}\quad
\UU(B_{-j}\otimes \omega_E(j))=
&B_{-j}\otimes \wedge^jU
}$$
where $U$ denotes the tautological sub bundle on the Grassmanian.
With this indexing, we will show that the maps $(\Psi_\F)_{i,j}$
and $(\Psi_\F)_{j,i}$ are dual up to a certain sign. 

By EFS Theorem 4.1 we may identify $A_{c+i}$ with
$\H^{k-i}(\F(i-k-1))$
and $B_j$ with $\H^j(\F(-j))$. As we have assumed that
$\F\iso \Ext^c(\F, \omega_{\P^n}(k+1))$, we have
$$
B_j^*=\H^j(\F(-j))^*
=\H^{j}(\Ext^c(\F, \omega_{\P^n}(k+1))(-j))^*
=\H^{n-j}(\F(k+1-j))
=A_j
$$
by Serre duality. With this identification it suffices
to check the 
signs in the maps $\varphi_P$ rather than in the maps
$\varphi_{\bullet,\bullet}$.

The linear complex
$\cH om(L,\omega_{\PP^n})(k+1)[c]$, is 
a linear free monad for the dual sheaf
$D(\F)\iso \F$. By the uniqueness and functoriality
of linear monads, the isomorphism $\sigma$ induces
an isomorphism
$L\iso \cH om(L,\omega_{\PP^n})(k+1)[c]$.

To simplify notation we set
$\check L^i = D(L^i)= \cH om(L^i,\omega_{\PP^n})(k+1).$
We follow standard sign conventions (see for example Iverson [1986])
and define the dual complex 
$\check L = \cH om(L,\omega_{\PP^n})(k+1)$ to 
have differentials $(-1)^i\check\alpha_i$. Shifting the 
complex $c$ steps also introduces the sign $(-1)^c$. 
Thus the isomorphism $L\to \check L[c]$ consists of 
a sequence of isomorphisms $\sigma_j: L^j\to \check L^{-c-j}$
as in the following diagram:
$$\diagram
 \rTo & L^{-c-i}& \rTo^{\alpha_{-c-i}} & \ldots & 
 \rTo^{\alpha_{-c-1}} & L^{-c}& \rTo & \ldots &
 \rTo^{\alpha_{j-1}} & L^{j}& \rTo  \cr
  & \dTo^{\sigma_{-c-i}}&  & & 
  & \dTo^{\sigma_{-c}}&  & &
  &\dTo^{\sigma_{j}}&   \cr
 \rTo & \check L^{i}& \rTo^{(-1)^{c+i-1}\check \alpha_{i-1}} & \ldots & 
 \rTo^{(-1)^c \check \alpha_{0}} & \check L^{0}& \rTo & \ldots &
 \rTo^{(-1)^{-j} \check \alpha_{-c-j}} & \check L^{-c-j}& \rTo  \cr
\enddiagram
$$
{}From the diagram we see that
$$\leqno{(*)}
\eqalign{
\varphi_{-c-i,j} =& \sigma_j \alpha_{j-1} \tensor \ldots \tensor \alpha_{-c-i} 
\cr
=&(-1)^s \check \alpha_{-c-j} \tensor \ldots 
\tensor \check \alpha_{i-1}\sigma_{-c-i}
}
$$
with $s=c(c+i+j)+{c+j+1 \choose 2}+{i \choose 2}$,
where the $c(c+i+j)$ comes from the shift, and the rest
is the contribution of the signs in the duality, separating
the parts with positive and negative indices.

We next prove that the map $\sigma_{-c-i}$ is, up to a sign
we shall identify,
the dual of $\sigma_i$. By the uniqueness and functoriality
of linear monads and the $\epsilon$ symmetry
of $\sigma$ the induced map of complexes $\sigma': L\to \check L[c]$
factors as the composite $\sigma' = \epsilon D(\sigma') \iota'$
where $\iota'$ is 
the canonical morphism of complexes
$$
\iota': L \rTo \cH om(\cH om(L,\omega_{\PP^n})[c],\omega_{\PP^n})[c].
$$
The components of $\iota'$ are given by
$$
\iota'_\ell: L_\ell \rTo^{(-1)^{(c+1)(c+\ell)} \iota } \check {\check L}_\ell,
$$
where $\iota$ denotes the canonical morphism $M \rTo \check {\check M}$
 of sheaves, c.f. Iverson [1981], p.73. Thus
$\sigma_{-c-i}= 
\epsilon (-1)^{(c+1)i} \check \sigma_i \iota$.

Combining this equation with $(*)$ we get
$$
\eqalign{
\varphi_{-c-i,j}
& =\epsilon (-1)^{s+(c+1)i} \check \alpha_{-c-j} \tensor \ldots 
\tensor \check \alpha_{i-1} \check \sigma_i \iota \cr
& =\epsilon (-1)^{s+t} {\rm transpose}(\sigma_i \alpha_{i-1} \tensor
\ldots \tensor \alpha_{-c-j}) \cr
}$$
with $t=(c+1)i+{c+i+j \choose 2}$ since in the transpose matrix the
tensor factors occur in the opposite order, and this tensor lies in
$\wedge^{c+i+j} W$.

Now
$$
\eqalign{
s+t&=c(c+i+j)+{c+j+1 \choose 2}+{i \choose 2}+(c+1)i+{c+i+j \choose 2} \cr
&= c^2+c(i+j)+{(c+j)^2+(c+j) \over 2}+{i^2-i \over 2}+
ci+i+ {(c+i+j)^2-(c+i+j) \over 2} \cr
&\equiv (c+1)(i+j)+ij \hbox{ mod } 2.
}$$
By \ref{gr signs}, we see that all the diagonal
blocks $(\Psi_\F)_{i,i}=\UU(\varphi_{-c-i,i})$
will be $\epsilon$ symmetric. 
Because
$(c+1)(i+j)+ij+k(i+j)+ij \equiv v(i+j)$
we can multiply the block matrix $\Psi_\F$
by the diagonal matrix of signs
$\Delta=\oplus_j(-1)^{vj}Id_{A_j}$
where $Id_{A_j}$ is the identity map
on $A_j$, to get
a map which is $\epsilon$ symmetric; that is,
setting 
$\Psi_\F'=
\Psi_\F\Delta
$
we will have
$
\epsilon \cH om((\Psi_\F'),\O_\GG(-1))
=\Psi_\F'.
$
\Box

\section{curves} Curves

Resultants of binary forms were the starting point for this subject
(Leibniz [1692],  B\'ezout [1779], 
Sylvester [1840], [1842], see Kline [1972] for some historical remarks), and they
correspond to the simplest cases of Chow forms of curves.
We begin by explaining how they fit into our theory.

\medskip

\example{binary}{\bf Binary Forms} 
Consider the {\bf rational normal curve} 
$\PP^1 \hookrightarrow \PP^d, \, (s:t) \mapsto (s^d,s^{d-1}t,\ldots, t^d)$
we use 
$[i,j]$ for the $ij^{th}$ Pl\"ucker coordinate of 
$\GG_2=\GG(2,\H^0(\PP^d,\cO(1)))$
with respect to the given basis. The following determinantal
formula can be deduced directly by computation from the 
classic B\'ezout formula for the resultant, or, since
the rational normal curve is a linear determinantal variety,
it could be deduced from
\ref{chow of Ulrich}.

{}From our point of view the most direct method is the  computation of a map 
in a Tate resolution.

\proposition{P1}
The Chow form of the rational normal curve of degree $d$
is the determinant of the $d \times d$ symmetric matrix
$A=(a_{ij})$ with
$$
a_{ij} = \sum_{p < min(i,j) \atop p+q=i+j-1} [p,q].
$$

\proof Consider $\cL=\cO_{\PP^1}(-1)$. 
With respect to $H=\cO_{\PP^1}(d)$ the betti numbers of the Tate
resolution of $\cL$ are
$$
\matrix{ 3d & 2d & d & - & - & - \cr
          - &  - & - & d & 2d & 3d \cr }
$$
Let $y_0,\ldots,y_d$ denote the dual basis
in $V$ to the monomial basis in $W=\H^0\O_{\P^1}(d)$.

The $d \times 2d$ matrix comes from 
multiplication $\H^0(\PP^1,\cO(d-1))\times \H^0(\PP^1,\cO(d)) \to
\H^0(\PP^1,\cO(2d-1))$, hence is given by the Sylvester type matrix
$$
B=(b_{kl})=(y_{k-l})=\hbox{ transpose }\pmatrix{ y_0 & y_1 & \ldots  & y_d & 0 & \ldots & 0 \cr
          0 & y_0 & \ldots & y_{d-1} & y_d &\ldots & 0 \cr
          \vdots&\ddots&\ddots& &\ddots &\ddots &\vdots \cr
          0 & \ldots & 0 & y_0 & \ldots & y_{d-1} & y_d \cr}.
$$
To prove the formula we must show that the kernel of this matrix is the image
of a matrix $A=(a'_{i,j})$ with
$a'_{i,j}=\sum_{p < min(i,j) \atop p+q=i+j-1} y_p\wedge y_q.$

The equation $B\cdot A=0$ holds since a 
term $y_{k-l} \wedge y_p \wedge y_q $ arising in the product
$b_{kl}a'_{lj}$ is cancelled either by a term $y_p \wedge y_{k-l} \wedge y_q$
in the product $b_{k,k-p}a'_{k-p,j}$ or by a term $y_{q} \wedge y_{p}
\wedge y_{k-l}$ in $b_{k,k-q}a'_{k-q,j}$, in case $k-l < j$ or $j \le k-l$
respectively. 

Since the $d$ rows of $A$
are linearly independent, and we know that the kernel of $B$ is
 generated by $d$ independent elements of degree 2, we see that
the rows of $A$ generate the kernel of $B$ as required.
\Box

If we want to obtain the 
Sylvester resultant formula, we apply $\UU_2$ to the Tate
resolution shifted. The resulting complex 
$$ 
\oplus^d \cU \to \oplus^{2d} \cO, 
$$
written in Stiefel coordinates,
gives the classical Sylvester formula for two polynomials
$ f=f_0s^d+f_1s^{d-1}t+\ldots+f_dt^d$ and 
$g=g_0s^d+g_1s^{d-1}t\ldots+g_dt^d$ of equal degree.

We will generalize these formulas to arbitrary curves over
an algebraically closed field, and obtain partial results for
more general ground fields.

By a {\it curve} we will mean a purely one-dimensional 
scheme $X$, projective over $K$.
The theory of Ulrich sheaves on curves is significantly simpler
than the theory for higher-dimensional varieties because
it is essentially independent of the embedding. To state the
result, we say that a sheaf $\G$ on a curve $X$ has 
{\it no cohomology\/}
if $\H^0(\G)=\H^1(\G)=0$.

\theorem{arbitrary curves} If $X$ is a curve embedded in
$\P=\P^{n+1}$ with hyperplane divisor $H$, then a sheaf
$\F$ is an Ulrich sheaf for $X$ in $\P$ if and only if
$\F=\G(H)$ for some $\G$ with no cohomology.

\proof $\G(H)$ is 0-regular because $\H^1(\G(H)(-H))=\H^1(\G)=0$.
Similarly, $\cE xt^{n-1}_\P(\G(H),\O_\P(-n-1))$ is
2-regular because $\H^0(\G)=0$. (One can also see the
desired vanishing
directly from the Tate resolution: for example, the vanishing of
$\H^0(\G)$ implies that the free module
$T^0(\G)$ has no generators in degree 0; and it follows that
for $j<0$ the module $T^{-j}(\G)$ has no generators in
degree $-j$. But by Eisenbud-Fl\o ystad-Schreyer [2000, Thm.~7.1] the space of
generators of $T^{-j}(\G)$ in
degree $-j$ is $\H^0(\G(-jH)$.)\Box

To find sheaves with no cohomology
it suffices to look
for sheaves on a single component of the reduced scheme $X_{\rm red}$
or even on its normalization. Thus we are led to ask: Given a
nonsingular irreducible curve $X$ over an arbitrary field $K$,
what are the sheaves $\G$ over $X$ with no cohomology? Such a sheaf
$\G$ can have no torsion, so (since $X$ is nonsingular) 
$\G$ is automatically locally free. From the vanishing of the
cohomology we see that the Euler characteristic of $\G$ is 0,
so by Riemann-Roch the degree of $\G$ is $\rank(\G)\cdot (1-g)$,
where $g={\rm genus}(X)$. Over an algebraically closed field,
there are always line bundles of this type.
This generalizes the fact that the equation of any plane curve can be
written as the determinant of a linear matrix:

\proposition{ineff} A line bundle $L$ on a curve $X$ has no cohomology
 if and only if $\deg(L)={\rm genus}(X)-1$ and
$L$ has no sections. If $X$ contains infinitedly many $K-rational$
points,
then such line bundles exist on $X$, and thus the Chow form
of $X$, in any projective embedding, can be written as a
determinant of linear forms in the Pl\"ucker coordinates.

\proof The first statement is immediate from the Riemann-Roch theorem.
For the second, take $L=\O_X(p_1+\cdots+p_g-q)$, where the
$p_i$ and $q$ are general $K$-rational points.\Box

To arrive at explicit resultant formulas further work has to be done.
We have to compute the appropiate  differentials in the Tate complex
explicitedly.

\example{hyperelliptic resultant} {\bf Hyperelliptic resultant formulas} 
Consider a fixed polynomial $f=f_0+f_1t+\ldots+f_{2g+2}t^{2g+2}$ with no 
multiple roots.
To write explicit Stiefel and B\'ezout formulas for
the resultant of
two functions
$a(t)+b(t)\sqrt{f(t)}$ and $c(t)+d(t)\sqrt{f(t)}$ with $a,b,c,d  \in K[t]$  
we consider them as functions on the hyperellipic curve $C$ of genus g 
with function
field $K(t,\sqrt f )$.  Let  $k= \max\{\deg a, g+1+\deg b, \deg c,
g+1+\deg d \}$ and consider the embedding of $C$ given by
$t \mapsto (1:t:\ldots:t^k:\sqrt f:t\sqrt f:\ldots: t^{k-g-1}\sqrt f)$.
We want the Chow form of this embedding. By \ref{arbitrary curves}
and \ref{symmetry}, a formula as the determinant of a symmetric matrix
will arrize if choose as Ulrich sheaf $\L(H)$ with the line bundle $\L$ a 
``non vanishing theta characteristic" ---that is, a line bundle
$\L$ on $C$ such that $\L\otimes \L=\omega_C$, the
canonical bundle, and $\L$ has no cohomology. 
A non vanishing theta characteristic in turn corresponds to a 
factorization $f=f^{(1)}f^{(2)}$ of $f$ into two 
polynomials of degree $g+1$. All of our formulas will depend on the
choice of such factorization and we will obtain
 ${1\over 2}{2g+2\choose g+1}$  B\'ezout formulas.

Before we come to the  B\'ezout formulas we will prove a Stiefel
formula for the resultant that is highly parallel to the
Sylvester formula for the ordinary resultant. We will then
deduce a B\'ezout formula in a way that is analogous to our
proof of \ref{P1}.
Let
$$
syl(k,r) = \hbox{ transpose } \pmatrix{ r_0 & r_1 & \ldots  & r_k & 0 & \ldots & 0 \cr
          0 & r_0 & \ldots & r_{k-1} & r_k &\ldots & 0 \cr
          \vdots&\ddots&\ddots& &\ddots &\ddots &\vdots \cr
          0 & \ldots & 0 & r_0 & \ldots & r_{k-1} & r_k \cr}
$$
be the $2k\times k$  ``Sylvester block'' of a polynomial $r$ of
degree $k$.

\proposition{Hyperelliptic Sylvester formula}
With notation as above, two functions
$a+b\sqrt f$ and $c+d\sqrt f$ with $a,b,c,d \in K[t]$ 
have a common zero if and only if the determinant of the $4k \times 4k$ matrix
$$ \pmatrix{
        syl(k,a) & syl(k,bf^{(2)}) & syl(k,c) & syl(k,df^{(2)}) \cr  
        syl(k,bf^{(1)}) & syl(k,a) & syl(k,df^{(1)}) & syl(k,c) \cr }$$
vanishes. 

\proof Let $\pi: C \to \PP^1$ denote the double
cover corresponding to the inclusion $K(t)\subset K(C)=K(t)[\sqrt f]$.
We consider the embedding of $C$ as a curve of degree $2k$ 
in projective space $\PP^{2k+1-g}$
corresponding to the line bundle $\O_C(H)=\pi^*(\O_{\P^1}(k))$. 
The space of global sections of $\O_C(H)$ has basis corresponding
to the functions
$1,t,\ldots,t^k,\sqrt f,t\sqrt f,\ldots, t^{k-g-1}\sqrt f$,
so the Chow form of $C$ in this embedding is the resultant we seek.
We write $e_0,\ldots,e_k,e_{k+1},\ldots,e_{2k-g}\in V=\H^0(\O_C(H))^*$
for the dual basis.

Every line bundle $\L$ on $C$ can be described as a rank 2 vector bundle
$\B=\pi_* \L$ on $\PP^1$ together with an action $\B \rTo^y \B(g+1)$ satisfying
$y^2= f id_\B$. For example $\pi_* \O_C = \O \oplus \O(-g-1)$ with the action
defined by $y=\pmatrix{ 0 & f \cr 1 & 0 \cr}$. The bundle
$\B=\O(-1)\oplus\O(-1)$ with the action of 
$\pmatrix{0& f^{(1)} \cr f^{(2)} & 0\cr}$ corresponds to a non vanishing 
theta characteristic $\F$ on $C$. In particular, $\F$ is a line bundle of
degree $g-1$ with no cohomology.
See Buchweitz and Schreyer [2002] for a detailed exposition.
The Stiefel
formula above is obtained by applying the functor $\UU$ to the line bundle
$\F(2H)$.

The global sections of $\F(H)$ has a basis corresponding to
the functions
$$
\sqrt{f^{(1)}},t\sqrt{f^{(1)}},\ldots,t^{k-1}\sqrt{f^{(1)}},
\sqrt{f^{(2)}},t\sqrt{f^{(2)}},\ldots,t^{k-1}\sqrt{f^{(2)}},
$$
while  $\H^0(\F(2H))$ has a basis corresponding to
$$
\sqrt{f^{(1)}},t\sqrt{f^{(1)}},\ldots,t^{2k-1}\sqrt{f^{(1)}},
\sqrt{f^{(2)}},t\sqrt{f^{(2)}},\ldots,t^{2k-1}\sqrt{f^{(2)}}.
$$
Thus the map
$$
\Hom(E,\H^0(\F(H))) \to \Hom(E,\H^0(\F(2H)))
$$
in the Tate resolution is given 
by the $4k\times 2k$ matrix over the exterior algebra
$$
B=\pmatrix{\scriptstyle syl(k,e_0+e_1t+\ldots+e_kt^k) &\scriptstyle
syl(k,(e_{k+1}+\ldots+ e_{2k-g}t^{k-g-1})f^{(2)}) \cr
\scriptstyle syl(k,(e_{k+1}+\ldots+
e_{2k-g}t^{k-g-1})f^{(1)}) &\scriptstyle syl(k,e_0+e_1t+\ldots+e_kt^k) }.
$$
The desired Sylvester formula follows by interpreting the induced map
$$
\H^0(\F(H)) \tensor U \rTo \H^0(\F(2H)) \tensor \O_\GG
$$
in terms of Stiefel coordinates.
\Box

We now use these constructions as in \ref{P1}
to derive {\bf hyperelliptic B\'ezout
formulas}. It suffices to compute the kernel of 
the map $B$. By \ref{first main} this will be a
$2k \times 2k$ matrix with entries in $\Lambda^2 V$.
Because $\F$ is a theta characteristic
 \ref{symmetry} show that the kernel will be represented
in suitable bases by a symmetric matrix.

The final formula may be written in terms of the $2\times 2$
minors of the $2\times (2k+1-g)$ matrix 
$$
\pmatrix{ a_0 &\ldots & a_k & b_0 &\ldots & b_{k-g-1}  \cr
c_0 &\ldots & c_k & d_0 &\ldots & d_{k-g-1} \cr}.
$$ 
However we will work with the larger $2 \times 3(k+1)$ matrix
$$
\pmatrix{ a_0 &\ldots & a_k & (bf^{(1)})_0 &\ldots & (bf^{(1)})_k & 
(bf^{(2)})_0 &\ldots & (bf^{(2)})_k \cr
c_0 &\ldots & c_k & (df^{(1)})_0 &\ldots & (df^{(1)})_k & 
(df^{(2)})_0 &\ldots & (df^{(2)})_k \cr}
$$
whose minors are linear combinations of those of the matrix above
with coefficients, 
which depend on the coefficients of $f^{(1)}$ and $f^{(2)}$.

If $0 \le p,q \le k$ then we denote by $[p,q]$ the
 minor formed by the columns with indices $p$ and $q$. 
We write
$p^{(1)}$ for the column with index $p+(k+1)$, and $q^{(2)}$
the column with index $q+2(k+1)$. Thus brackets like
$[p^{(1)},q]$ and $[p^{(1)},q^{(2)}]$ represent $2\times 2$ minors of
the large matrix. 

Consider the $k \times k$ matrices $A^{11},\ldots,A^{22}$ defined by
$$
 \eqalign{ A^{11}_{i,j}& =
\sum_{0\le p<q \le k \atop {p < \min(i,j) \atop p+q=i+j-1}} [p^{(2)},q]+[p,q^{(2)}], \cr
A^{12}_{i,j}&=
\sum_{0\le p<q \le k \atop {p < \min(i,j) \atop p+q=i+j-1}} [p,q] \quad
+ \quad \sum_{0\le p,q \le k \atop {p < j \atop p+q=i+j-1}} [p^{(1)},q^{(2)}]
,\cr
A^{21}_{i,j}&=
\sum_{0\le p<q \le k \atop {p < \min(i,j) \atop p+q=i+j-1}} [p,q] \quad
+ \quad \sum_{0\le p,q \le k \atop {p < j \atop p+q=i+j-1}} [p^{(2)},q^{(1)}] ,\cr
 A^{22}_{i,j}&=
\sum_{0\le p<q \le k \atop {p < \min(i,j) \atop p+q=i+j-1}} [p^{(1)},q]+[p,q^{(1)}]. \cr}
$$ 

The matrix $A$ is actually symmetric. This becomes visible if we expand the
expressions into brackets of the smaller $2\times (2k+1-g)$ matrix.

\proposition{hyperelliptic bezout}
Suppose $k \le 12$. The functions $a+b\sqrt f$ and $c+d\sqrt f$ have a common zero
if and only if the determinant of the matrix
$$A=\pmatrix{A^{11} & A^{12} \cr A^{21} & A^{22} \cr }.$$
vanishes.

The formula should certainly hold for any $k$; but as noted in the 
proof we have performed the necessary computations only up to $k=12$.

Notice that in case $b=d=0$ the matrix reduces to twice the 
Bezout matrix 
for binary forms of degree $k$. This fits with the fact that 
two functions on $\P^1$ with a common zero have two common 
zeroes when  pulled back to $C$. 

\proof As in the proof of \ref{P1} it suffices to check that
$B\cdot A = 0$, when we regard $A$ as a matrix over the exterior algebra,
because the linear independence of the columns of $A$ is visible from
the specialization to the case of binary forms  $b=d=0$.
For
each specific value of $g$ and $k$ this can checked by Computer algebra,
and we did this for all cases $1\le g+1 \le k \le 12$. \Box

As a concrete application of \ref{hyperelliptic bezout} we do
the case of an elliptic curve over the complex numbers. 

\noindent
\example{resultant of doubly periodic functions}{\bf Resultant
of doubly periodic functions}
Consider an elliptic curve $C=\CC/\Gamma$ and the corresponding
Weierstrass $\wp$-function, with functional equation
$$
\wp'(z)= 4\wp^3(z)-g_2\wp(z)-g_3=4(\wp(z)-\rho_1)(\wp(z)-\rho_2)(\wp(z)-\rho_3)
$$
where the ${\rho_j}'s$ are the values of $\wp$ at the half periods.
Two doubly periodic functions
$$
f(z)= a_0+a_1 \wp(z) + a_2 \wp^2(z)+ b_0 \wp'(z)/2
$$
and 
$$
g(z)= c_0+c_1 \wp(z) + c_2 \wp^2(z)+ d_0 \wp'(z)/2$$
have a common zero iff the determinant of
$$
\pmatrix{ 
-\rho_1 \rho_2 [13]-(\rho_1+\rho_2) [03]&
      -\rho_1 \rho_2 [23]+[03]&
      [01]&
      [02]\cr
      -\rho_1 \rho_2 [23]+[03]&
      (\rho_1+\rho_2) [23]+[13]&
      [02]&
      [12]\cr
      [01]&
      [02]&
      \rho_3 [13]+[03]&
      \rho_3 [23]\cr
      [02]&
      [12]&
      \rho_3 [23]&
      -[23]\cr
      }$$
vanishes, where the bracket
$
[ij]
$ 
denotes the minor made from the $i^\th$ and $j^\th$
columns of the matrix
$$
\pmatrix{
 a_0 & a_1 & a_2 & b_0 \cr  
c_0 & c_1 & c_2 & d_0
}.
$$
This formula follows from \ref{hyperelliptic bezout},
with one of the roots of $f$ at infinity, and with the 
factorization given by $f^{(1)}= (\wp(z)-\rho_1)(\wp(z)-\rho_2)$.

\medskip

Returning to our general discussion, we may ask whether it
is possible to give a B\'ezout formula for the Chow
form of a curve over a field $K$ even if the curve does not have enough
$K$-rational
points to apply \ref{arbitrary curves}. In this case the curve
may have no rank 1
Ulrich sheaf, as happens, for example, for a conic without real points in 
$\P^2_{\bf R}$. However, it may be that there are always rank 2 Ulrich sheaves.
For example, assuming that $X$ has genus 0, 
The  structure sheaf $\O_X$ and the canonical bundle $\omega_X$
 are defined over $K$, and there is a unique
extension 
$$
\eta:\quad 0\to \omega_X\to \cE \to \O_X\to 0
$$
corresponding to a nonzero element $\eta\in \H^1(\omega_X^{-1}) = K$.
Over an algebraic closure of $K$ the bundle $\cE$ splits as
$\O_{\P^1}(-1)\oplus \O_{\P^1}(-1)$ (the sequence above
is the Koszul complex) and thus $\cE$ has no cohomology.

The main theorem of Brennan, Herzog, and Ulrich [1987] generalizes
this example and 
says that if $K$ is algebraically closed and 
$X$ is a 1-dimensional arithmetically Cohen-Macaulay subscheme
of $\P$ then there exists a rank 2 sheaf $\F$ with no cohomology, 
which in addition
satisfies $\F\iso\cH om(\F, \omega_X)$.
(Their statement does not include the separability hypothesis
below; but they apply a result
of Eisenbud [1988] which is proved only in the
algebraically closed case. We do not see how to extend their
proof beyond the separable case, as below.)
A variation
on their proof allows one to drop the ``arithmetically Cohen-Macaulay"
hypothesis. Here is a geometric version of the argument, 
developed in conversation with Joe Harris.

\proposition{inf field case} Let $X$ be a projective
curve, separable
over the field $K$.
If $K$ is infinite
then $X$ has a coherent sheaf $\cE$ with no cohomology which
is a rank 2 vector bundle over
the normalization of $X_{\rm red}$, 
and satisfies $\Hom(\cE, \omega_X)=\cE$.

\proof Let $\pi: C\to X_0$ be the normalization.
It is enough to find a rank 2 vector bundle without cohomology on $C$
with $\Hom(\cE, \omega_C)$, since we have
$\Hom_X(\cE, \omega_X)
=\Hom_C(\cE,\Hom(\O_C,\omega_X)=\Hom_C(\cE, \omega_C)$. 
Since we have dealt with the case of $\P^1$ above, we will
assume that the genus $g$ of $C$ is greater than 0.
Let $L$ be a line bundle on $C$ of strictly positive degree

Any extension class
$\eta
\in \ext^1(\omega_C\otimes L, L^{-1})$ gives rise to a short exact
sequence
$$
\eta:\quad 0\to L^{-1}\to \cE \to \omega_C\otimes L\to 0
$$
where $\cE$ is a vector bundle. For any such bundle 
$\wedge^2\cE=\omega_C$, whence $\Hom(\cE, \omega_C)=\cE$.

By Serre duality $\chi(\cE)=0$, so $\cE$ will be an Ulrich
sheaf as long as $\H^0(\cE)=0$. 
Since $\H^0(\L^{-1})=0$,  this condition is satisfied if and only
if the connecting homomorphism 
$$
\delta_\eta: \H^0(\omega_C\otimes L)\to \H^1(\L^{-1})=\H^0(\omega_C\otimes L)^*
$$
is an isomorphism. But 
$$
\eta\in \ext^1(\omega_C\otimes L, L^{-1})\iso \H^1(L^{-2}\otimes\omega_C^{-1})
\iso \H^0(L^2\otimes\omega_C^2)^*,
$$
and $\delta_\eta$ is induced by the multiplication pairing
$$
\H^0(L\otimes\omega_C)\otimes \H^0(L\otimes\omega_C)\rTo^m
\H^0(L^2\otimes\omega_C^2)
$$
in the sense $\eta$ goes to $\delta_\eta$ under the composite
$$\eqalign{
\H^0(L^2\otimes\omega_C^2)^*\rTo^{m^*}
&\H^0(L\otimes\omega_C)^*\otimes\H^0(L\otimes\omega_C)^*\cr
&\iso \H^0(L\otimes\omega_C)^*\otimes\H^1(L^{-1})\cr
&\iso \hom(\H^0(L\otimes\omega_C),\ \H^1(L^{-1})).
}$$
Now the ring $R=\oplus_d\H^0(L^d\otimes\omega_C^d)$
is an integral domain, by separability it splits
into a product of integral domains over the algebraic
closure of $K$. It follows that the 
multiplication pairing is a direct sum of 1-generic
pairings in the sense of
Eisenbud [1988]. The results of that paper show that
$\delta_\eta$ is an isomorphism unless $\eta$ lies in
a certain proper hypersurface in $\H^0(L^2\otimes\omega_C^2)$.
If $K$ is infinite then this hypersurface cannot contain
all the $K$ rational points of this vector space.\Box

\example{conic with no points} {\bf A Conic without a point} 
The conic $C \subset \PP^2$ defined by $x^2+y^2+z^2=0$ has no line bundle
of degree $-1$ defined over $\RR$. However there are rank 2 Ulrich sheaves.
The cokernel 
$$\F=\coker(\cO_{\PP^2}^4(-2) \rTo^M \cO_{\PP^2}^4(-1))$$
given by the matrix
$$M= \pmatrix{0 & x & y & z \cr
              -x & 0 & z & -y \cr
              -y & -z & 0 & x \cr
              -z & y & -x & 0 \cr}
$$
is a rank 2 sheaf on $C$ with no cohomology. An explicit 
formula can be derived from Pfaffian B\'ezout formula for the
resultant of 3 quadratic forms in 3 variables given in the
introduction, by specializing one of the three quadratic forms
to $x^2+y^2+z^2$ and eliminating unnecessary variables.

\section{Projective spaces} Veronese embeddings and Resultant formulas

Consider the d-uple embedding 
$$\PP^k \hookrightarrow \PP^N$$ with $N = {d+k \choose k}-1$. The Chow
form is the resultant of $k+1$ homogeneous forms of degree $d$ in
$k+1$ variables, hence is of particular interest.  To find
determinantal or Pfaffian formulas for powers of such Chow forms, we
need to look for vector bundles on $\PP^k$ that become Ulrich sheaves
on the $d$-uple embedding; Stiefel formulas come from weakly Ulrich
sheaves.  By an argument shown to us by Jerzy Weyman, even Ulrich
sheaves always exist! By way of comparison, the classical search for
B\'ezout or Stiefel formulas was essentially a search for line bundles
on $\PP^k$ that become Ulrich or weakly Ulrich on the $d$-uple
embedding. Weakly Ulrich line bundles exist (and were found
classically) if and only if $k\leq 4$ or $k=5,\; d\leq 3$ (Ulrich line
bundles never exist except when $k\leq 2$ or $d=1$.) We get a few more
Stiefel formulas for the resultants themselves (and not just powers)
from the Horrocks-Mumford bundle in the case $k=5,\; d=4,6$ or 8.

It turns out that the cohomology of a sheaf that becomes an Ulrich sheaf
on the $d$-uple embedding is
determined by the rank of the sheaf
alone, and the same idea works for the $d$-uple embedding
of any variety:

\theorem{d-uple} Let $\iota: \P^m\hookrightarrow \P^n$ 
be the $d$-uple embedding. Suppose $\F$ is a sheaf
of dimension $k$ on $\P^m$.
The sheaf $\iota_*\F$ is an Ulrich sheaf on $\P^n$ if and only
if
$$
h^i(\F(e))\neq 0 \Leftrightarrow
\cases{
i=0,\quad -d<e \cr
0<i<k,\quad -(i+1)d<e<-id\cr
i=k,\quad e<-kd}.
$$
In particular, $\F$ then
has natural cohomology as sheaf on $\P^m$. Thus all the 
$h^i(\F(e))$ are determined by the formula
$$
\chi(\F(e))=h^0(\F){{e\over d}+k\choose k}.
$$
If $\F$ is a vector bundle of rank $r$ on $\P^m$, then
we can rewrite this formula as
$
\chi(\F(e))=
{r \over m!}(e+d)\cdots(e+md)=
({r \over m!}e^m)+\cdots+rd.$

The vanishing and non-vanishing results in the first part of 
\ref{d-uple}
have a very simple interpretation in terms
of the  betti diagram of the 
Tate resolution of $\F$: they say that the nonzero terms
form a sequence of non-overlapping strands
and that all of the strands 
representing intermediate 
cohomology have length
precisely $d-1$. The formulas in the second part then give the 
values of the nonzero terms. For example, if $\F$ is 
a rank 2 vector bundle on $\P^2$ which is an Ulrich sheaf
for the $d$-uple embedding, 
\ref{d-uple}
says precisely that the Tate resolution of $\F$,
considered as a sheaf on $\P^2$, has betti diagram
$$
\matrix{
\scriptstyle\cdots&\scriptstyle 2(d+2)
&\scriptstyle 1(d+1)&\scriptstyle0&\scriptstyle0 
&\scriptstyle 0 
&\scriptstyle 0 &\scriptstyle0 &\scriptstyle0&\scriptstyle0
&\scriptstyle\cdots&\scriptstyle0&\scriptstyle\cdots\cr
\scriptstyle\cdots&\scriptstyle0&\scriptstyle0&\scriptstyle1(d-1)
&\scriptstyle2(d-2)
&\scriptstyle\cdots&\scriptstyle(d-2)2&\scriptstyle(d-1)1&\scriptstyle0
&\scriptstyle0&\scriptstyle\cdots&\scriptstyle0&\scriptstyle\cdots\cr
\scriptstyle\cdots&\scriptstyle0&\scriptstyle0&\scriptstyle0     
&\scriptstyle0   &\scriptstyle0     &\scriptstyle0     
&\scriptstyle0     &\scriptstyle1(d+1)&\scriptstyle2(d+2)
&\scriptstyle\cdots&\scriptstyle d(2d)&\scriptstyle\cdots&\cr
}
$$
where the zeroth term is the one occuring at the far right,
(so that for example $h^0(\F)=2d^2$). Further examples are
given in the discussion of sheaves on $\P^3$ below.

To prove that the cohomology vanishes as we claim, we will repeatedly
use the following elementary result, which is an easy case of
Eisenbud-Fl\o ystad-Schreyer [2000 Lemma 7.4]:

\lemma{generalized Mumford} Suppose $\G$ is a sheaf on $\P^k$.
\item{a)} If \/ $\H^{i+j}(\G(-1-j))=0$ for all $j\geq 0$ then $\H^i(\G)=0$.
\item{b)} If \/ $\H^{i-j}(\G(1+j))=0$ for all $j\geq 0$ then $\H^i(\G)=0$.

Note that the case $i=1$ in part a)
is Mumford's result showing that a
$(-1)$-regular sheaf is 0-regular. For the reader's convenience
we give a quick proof.
\medskip

\noindent{\sl Proof of \ref{generalized Mumford}.\/} a): 
Translating the condition in a) to a condition on the Tate
resolution $T^{\bullet}(\G)$
over the exterior algebra $E$, we see that the free summand
$\H^i(\G)\otimes \omega_E$ in $T^0(\G)$ maps injectively into
$T^1(\G)$. Since $T^{\bullet}(\G)$ is a minimal complex
and $E$ is Artinian, this is impossible. 

Part b) follows by applying the same argument to the dual of
the Tate resolution.\Box
\medskip

\noindent{\sl Proof of \ref{d-uple}.\/}
We begin by showing by induction on $i$ 
that for $i<k$ we have $\H^i(\F(e))=0$
if $e\leq -(i+1)d$. Since $\F$ becomes an Ulrich sheaf
under the $d$-uple embedding
we have $\H^0(\F(-d))=0$, and it follows that 
$\H^i(\F(e))=0$ for $e\leq -d$, which is the case $i=0$.
For $i>0$ we proceed by descending induction on $e$. Again
since $\F$ becomes Ulrich on the $d$-uple embedding we have
$\H^i(\F(-(i+1)d))=0$, the initial case.
Assuming that 
$\H^i(\F(e))=0$ for some $e<-(i+1)d$,
the induction on $i$ gives
the hypothesis to apply part b) of
\ref{generalized Mumford}
to show that $\H^i(\F(-e-1))=0$.

Similarly, 
$\H^i(\F(e))=0$
for $i>0$ and $e\geq -id$ follows by induction and part b) 
of \ref{generalized Mumford}. The nonvanishing of the remaining
cohomology follows, since otherwise the Tate resolution for $\F$
would contain terms equal to zero. 

We next prove the formulas for $\chi(\F(e))$.
If $\iota_*\F$ is an Ulrich sheaf,
then  \ref{coho of Ulrich}
shows that
$\chi(\F(dt))=h^0(\F){k+t\choose k}$. Since
$\chi(\F(t))$ is a polynomial, it is determined by
this relation, yielding the first formula.

If in addition $\F$ is a bundle of rank $r$ on $\P^k$,
then part c) of \ref{Ulrich char}
shows that $\h^0(\F)=\deg\iota_*(\F)$,
which is $r$ times the degree of the $d$-uple embedding of
$\P^k$, that is, $\h^0(\F)=rd^k$. Substituting this in the
first formula we get the last formulas. (One could also
argue directly from the fact that the last formula must be
a polynomial of degree $k$ which vanishes at $-nd$ for
$n=1,\dots,k$.)
\Box

\corollary{divisibility} Suppose there exists a rank $r$ sheaf on
$\P^k$ which is an Ulrich sheaf for the $d$-uple embedding.
If a prime $p$ divides $d$ and $p^t$ divides $k!$,
then $p^t$ divides $r$. For example, 
any Ulrich sheaf on the $k!$-uple embedding of $\PP^k$
has rank a multiple of $k!$.

\proof In \ref{d-uple}, note that $\chi(\F(1))$
is an integer.\Box

The general problem of finding (weakly) Ulrich sheaves for the
Veronese embeddings of a given variety can be reduced to problem for
projective spaces by using a finite projection map (that is, a finite
map such that $\pi^*(\O_{\P^k}(1) = \O_X(1)$---such things always
exist by ``Noether normalization'') onto a projective space. 
This result and the following Corollary were 
inspired by the proof of the existence of rank
4 Ulrich sheaves on the 4-uple embedding of $\P^3$ given
by Douglas Hanes in his thesis [1999].

\proposition{general veronese} Let $X\subset \P^n$ be
a purely $k$-dimensional scheme, and let $\F$ be an
Ulrich sheaf whose support is $X$. Suppose that
$\pi: X\to \P^k$ is a finite projection. If 
$\E$ is a sheaf on projective space that is (weakly)
Ulrich for the $d$-uple embedding of projective space,
then $\F\otimes \pi^*\E$ is (weakly) Ulrich for the 
$d$-uple embedding of $X$.

\proof Since the cohomology of  $\pi_*\F(n)$ is the same as
the cohomology of $\F(n)$, we see from the cohomological 
characterization of Ulrich sheaves that $\pi_*\F$ is a trivial
bundle $\O_{\P^k}^t$ on $\P^k$. Since
$$
\H^q(\F\otimes \pi^*\E(d))=\H^q(\pi_*\F\otimes \E(d)),
$$
this group vanishes for exactly the same values of $q,d$ as
does $\H^q(\E(d))$, and this determines the weakly Ulrich and 
Ulrich properties.
\Box

If we apply \ref{general veronese} 
in the case where $X\iso \P^k$, embedded by the $e$-uple
embedding, we get a weak converse to \ref{divisibility}.

\corollary{semigroup Ulrich} If $\P^k$ has Ulrich sheaves of
ranks $a$ and $b$ on its
$d$-uple and $e$-uple embeddings respectively, then it has an Ulrich
sheaf of rank $ab$ on its $de$-uple embedding.\Box

If our ground field $K$ has characteric zero then 
any indecomposable homogenous bundle
on $\PP^n$ can be obtained by applying a Schur functor $S_\lambda$ to
the universal rank n quotient bundle $Q=\coker\O_{\PP^n}(-1) \to \O_{\PP^n}^{n+1}$ 
of $\PP^n$ (the tangent bundle tensor $\O_{\PP^n}(-1)$).
Here $\lambda=(\lambda_1,\ldots,\lambda_n)$ is a partition into at most
$n$ parts. Note $(S_\lambda Q)(1)=S_{\lambda+(1,1,\ldots,1)} Q$ 
and $\H^0(S_\lambda Q) = S_\lambda V$ with
$V = H^0(\cO(1))^*$.
Thus up to twist we may assume that $\lambda_n=0$. The theorem below implies that
$S_\lambda Q$ has Castelnuovo-Mumford regularity precisely zero iff $\lambda_n=0$.
For our purposes it is convenient to visualize the partition as a Ferrers diagram
whose row lengths are given by the $\lambda_i$, as follows:
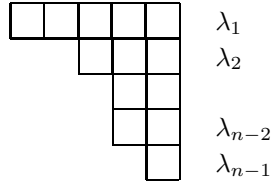
\centerline{
\hbox { \vbox { \vskip 5pt
        \hrule width 64pt
        \vskip 0pt
        \hbox {
                \kern -3.5pt
                \hbox to 13pt{\vrule height10pt depth3pt$\ $ \vrule height10pt depth3pt}
                \kern -3.5pt
                \hbox to 13pt{\vrule height10pt depth3pt$\ $ \vrule height10pt depth3pt}
                \kern -3.5pt
                \hbox to 13pt{\vrule height10pt depth3pt$\ $ \vrule height10pt depth3pt}
                \kern -3.5pt
                \hbox to 13pt{\vrule height10pt depth3pt$\ $ \vrule height10pt depth3pt}
                \kern -3.5pt
                \hbox to 13pt{\vrule height10pt depth3pt$\ $ \vrule height10pt depth3pt}  $\quad \lambda_{1}$
                }
        \vskip 0pt
        \hrule width 64pt
        \vskip 0pt
        \hbox { \hskip 26pt \kern -3.5pt
                \vbox{ \hbox{
                        \kern -3.5pt
                        \hbox to 13pt{\vrule height10pt depth3pt$\ $ \vrule height10pt depth3pt}
                        \kern -3.5pt
                        \hbox to 13pt{\vrule height10pt depth3pt$\ $ \vrule height10pt depth3pt}
                        \kern -3.5pt
                        \hbox to 13pt{\vrule height10pt depth3pt$\ $ \vrule height10pt depth3pt}  $\quad \lambda_{2}$
                        }
                \vskip 0pt
                \hrule width 38pt
                \vskip 0pt
                \hbox { \hskip 13pt \kern -3.5pt
                      \vbox { \hbox { \kern -3.5pt 
                                   \hbox to 13pt{\vrule height10pt depth3pt$\ $ \vrule height10pt depth3pt}
                                   \kern -3.5pt
                                   \hbox to 13pt{\vrule height10pt depth3pt$\ $ \vrule height10pt depth3pt} 
                               }
                               \vskip 0pt
                               \hrule width 25pt
                               \vskip 0pt 
                               \hbox { \kern -3.5pt 
                                    \hbox to 13pt{\vrule height10pt depth3pt$\ $ \vrule height10pt depth3pt}
                                    \kern -3.5pt
                                    \hbox to 13pt{\vrule height10pt depth3pt$\ $ \vrule height10pt depth3pt} $\quad \lambda_{n-2}$
                                }
                               \vskip 0pt
                               \hrule width 25pt
                               \vskip 0pt
                               \hbox { \hskip 13pt  \kern -3.5pt
                                     \vbox{ \hbox{ \kern -3.5pt
                                      \hbox to 13pt{\vrule height10pt depth3pt$\ $ \vrule height10pt depth3pt} $\quad \lambda_{n-1}$}
                                      \vskip 0pt
                                      \hrule width 12pt
                                      \vskip 0pt
                                      }}}}}}}}}

\noindent
The following result was pointed out to us by J. Weyman.

\theorem{Tate of homogeneous bundles} Suppose that $K$ has characteristic zero. 
Let $\lambda=(\lambda_1,\ldots,\lambda_{n-1})$ be a partition and $Q$
the universal rank n quotient bundle on $\PP^n$. 
The Tate resolution of the homogeneous bundle $\F = S_\lambda Q$ has nonzero terms
only where there are $*$s in the following diagram, in which the Ferrers diagram has
shape $\lambda$ as above:
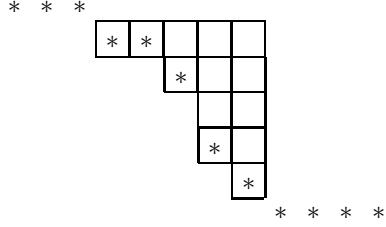
\vbox{ \vskip 5pt  \hbox { \hskip 80pt \vbox {
\vskip 5pt $ * $ \hskip 4pt $ *$ \hskip 4pt $ *$ \vskip 2pt
\hbox { \hskip 50pt \vbox{
        \hrule width 64pt
        \vskip 0pt
        \hbox {
                \kern -3.5pt
                \hbox to 13pt{\vrule height10pt depth3pt$\ *$ \vrule height10pt depth3pt}
                \kern -3.5pt
                \hbox to 13pt{\vrule height10pt depth3pt$\ *$ \vrule height10pt depth3pt}
                \kern -3.5pt
                \hbox to 13pt{\vrule height10pt depth3pt$\ $ \vrule height10pt depth3pt}
                \kern -3.5pt
                \hbox to 13pt{\vrule height10pt depth3pt$\ $ \vrule height10pt depth3pt}
                \kern -3.5pt
                \hbox to 13pt{\vrule height10pt depth3pt$\ $ \vrule height10pt depth3pt}
                }
        \vskip 0pt
        \hrule width 64pt
        \vskip 0pt
        \hbox { \hskip 26pt \kern -3.5pt
                \vbox{ \hbox{
                        \kern -3.5pt
                        \hbox to 13pt{\vrule height10pt depth3pt$\ *$ \vrule height10pt depth3pt}
                        \kern -3.5pt
                        \hbox to 13pt{\vrule height10pt depth3pt$\ $ \vrule height10pt depth3pt}
                        \kern -3.5pt
                        \hbox to 13pt{\vrule height10pt depth3pt$\ $ \vrule height10pt depth3pt}
                        }
                \vskip 0pt
                \hrule width 38pt
                \vskip 0pt
                \hbox { \hskip 13pt \kern -3.5pt
                      \vbox { \hbox { \kern -3.5pt 
                                   \hbox to 13pt{\vrule height10pt depth3pt$\ $ \vrule height10pt depth3pt}
                                   \kern -3.5pt
                                   \hbox to 13pt{\vrule height10pt depth3pt$\ $ \vrule height10pt depth3pt}
                               }
                               \vskip 0pt
                               \hrule width 25pt
                               \vskip 0pt 
                               \hbox { \kern -3.5pt 
                                    \hbox to 13pt{\vrule height10pt depth3pt$\ *$ \vrule height10pt depth3pt}
                                    \kern -3.5pt
                                    \hbox to 13pt{\vrule height10pt depth3pt$\ $ \vrule height10pt depth3pt}
                                }
                               \vskip 0pt
                               \hrule width 25pt
                               \vskip 0pt
                               \hbox { \hskip 13pt  \kern -3.5pt
                                     \vbox{ \hbox{ \kern -3.5pt
                                      \hbox to 13pt{\vrule height10pt depth3pt$\ *$ \vrule height10pt depth3pt}}
                                      \vskip 0pt
                                      \hrule width 12pt
                                      \vskip 3pt \hbox {\hskip 16pt $*$ \hskip 4pt $*$ \hskip 4pt $*$  \hskip 4pt $*$}
                                      }}}}}}}}}}}         
More precisely, for $1\le i \le n-1$ the cohomology group 
$\H^i((S_\lambda Q)(m))$ is nonzero
 if and only if  
$\lambda_{n-i+1} < -m-i \le \lambda_{n-i}$,
 $\H^0 S_\lambda Q (m)=0$ iff $m<0$ and 
$\H^n S_\lambda Q(m)=0$ iff $m\ge -n-\lambda_1-1$.

\proof The cohomology of a homogeneous bundle on the homogeneous space  
$\PP^n = \GL(n+1)/\pmatrix{ \GL(n) & * \cr 0 & \GL(1) }$ is determined by  
Bott's formula, see Jantzen [1981]. In particular 
$\H^iS_\lambda Q(m) \not= 0 \hbox{ at most one } i$ and 
$$\H^iS_\lambda Q(m) = 0 \hbox{ for all } i \Leftrightarrow 
-m \in \{\lambda_i+n+1-i \mid i=1,\ldots,n \}.$$
Thus the Hilbert polynomial $\chi S_\lambda Q (m)$ has precisely $n$ integral zeroes
and the Tate resolution ``steps down'' precisely at these $n$ values   
by \ref{generalized Mumford}. 
\Box

\corollary{homogeneous Ulrich} Suppose that $K$ has characteristic zero. 
The unique
indecomposable homogeneous bundle on $\PP^n$ that is an Ulrich sheaf 
for the $d$-uple embedding is
$S_\lambda Q$ with $(\lambda=((d-1)(n-1),(d-1)(n-2),\ldots,(d-1),0)$. 
It has rank $d^{n \choose 2}$.

\proof The first statement follows easily from the previous Theorem. 
The rank of $S_\lambda Q$ is given by the hook formula( see Stanley [1971]
or Fulton [1997])
$$\rank S_\lambda Q = \prod_{(i,j) \in \lambda} {{n+i-j} \over h(i,j)},$$
where $h(i,j)$ denotes the hook length of the $(i,j)^{th}$ box. The largest hook lenght
is $h(1,1)=(d-1)(n-1)+(n-2)=d(n-1)-1$. The denominators of the first row contribute with
$\prod_j h(1,j)= (d(n-1)-1)(d(n-1)-2)\cdot\ldots\cdot(d(n-1)-d+1)(d(n-2)-1)
\cdot\ldots  \cdot 1=
{[d(n-1)]! \over d^{n-1} (n-1)!}$. The numinators give 
${[d(n-1)]!\over (n-1)! }$. Thus the first row contributes with $d^{n-1}$ and the
total poduct yields $d^{{n-1 \choose 2}+n-1}=d^{n \choose 2}$ by induction.
\Box

\medskip \noindent
{\bf Chow forms from line bundles on projective spaces}
\medskip

All the classically known formulas (and no new ones) 
for the resultant of
$k+1$ forms of degreed $d$ in $k+1$ variables come from
applying these ideas to line bundles on projective spaces.
We get B\'ezout formulas in this way only
for binary forms of
any degree or linear forms in any number of variables
by \ref{homogeneous Ulrich}.

On the other hand
$\cL= \cO(j)$ on $\PP^k$ gives rise to a 2 term complex, and hence a
Stiefel formula for the Chow form of the d-uple image, iff
$$ \H^0 \cL(-H) =0 \, \hbox{ and } \, \H^k \cL(-(k-2)H) = 0, $$
equivalently, iff
$$ d-1 \ge j \ge -k +(k-2)d=(k-2)(d-1)-2.$$
\noindent
Thus the Chow forms of $\PP^1,\PP^2,\PP^3$ for arbitrary $d$-uple embeddings, on
$\PP^4$ for quadrics and cubics, and on $\PP^5$ for quadrics, can be written
as determinants of maps of vector bundles on the Grassmannian, or as
determinantal formulas in the Stiefel coordinates. This is
precisely the list of 
Gel'fand, Kapranov and Zelevinski, [1994] Chap. 13, Prop. 1.6,
and the formulas are the same. For instance in the case of three ternary
quadrics we have:

\example{quadrics on P2}
For the 2-uple embedding (quadrics) of
$\PP^2$ the line bundle $\O_{\P^2}(1)$ is weakly Ulrich,
and we see that the Chow form is the determinant
of a canonical map on the Grassmannian $\GG$
$$
\O_\GG(-1)^6\to \U \oplus \O_\GG(-1)^3.
$$
The map is easy to calculate, and
in Stiefel coordinates it has matrix
$$\pmatrix{
a_0 & b_0 & c_0 & {\scriptstyle[0,1,5]}           &0                  &{\scriptstyle [0,1,2] }\cr
a_1 & b_1 & c_1 &{\scriptstyle [0,3,5]}           &{\scriptstyle [0,3,4]
}         & {\scriptstyle [0,1,4]-[0,2,3] }\cr
a_2 & b_2 & c_2 &{\scriptstyle [0,4,5]-[1,2,5]} &{\scriptstyle [0,3,5] }         & {\scriptstyle  [0,1,5] }\cr          
a_3 & b_3 & c_3 & 0                   &{\scriptstyle [1,3,4]}          &{\scriptstyle [0,3,4] }\cr       
a_4 & b_4 & c_4 &{\scriptstyle [2,3,5]}           & {\scriptstyle[2,3,4]+[1,3,5]} &{\scriptstyle [0,3,5]} \cr        
a_5 & b_5 & c_5 &{\scriptstyle [2,4,5]}           &{\scriptstyle [2,3,5]}          &0       \cr      
}.$$
Thus the determinant of this matrix
is the resultant of three quadratic forms
$d=d_0x^2+d_1xy+d_2xz+d_3y^2+d_4yz+d_5z^2$ for $d=a,b,c$
with $(i,j,k)\th$ Pl\"ucker coordinate
$$[i,j,k]=det \pmatrix{
a_i & b_i & c_i \cr
a_j & b_j & c_j \cr
a_k & b_k & c_k \cr}.$$

\goodbreak
\medskip \noindent
{\bf Ulrich sheaves on $\PP^2$}
\medskip

To get new formulas for resultants, we replace line bundles
with vector bundles of higher rank. The Chow forms of these
bundles are the desired resultants raised to a power equal
to the rank of the bundle. But if the rank the bundle is 2, then its
natural symplectic structure allows us to
find a polynomial square root by taking a Pfaffian in place
of a determinant, so we get formulas for the resultant itself.

\proposition{rank 2 on P2} If $\alpha$ is
a $(d+1)\times (d-1)$ matrix of linear forms on $\P^2$ 
whose minors of order $d-1$ generate an ideal of codimension 3
(the generic value), then
$$
\coker \bigl(\O_{\P^2}(d-2)^{d-1}\rTo^\alpha \O_{\P^2}(d-1)^{d+1}\bigr)
$$
is an Ulrich sheaf on the $d$-uple embedding of $\P^2$.

For example, we may take
$$
\alpha=\pmatrix{
x_0 & x_1 & x_2 & 0 & \ldots & 0 \cr
0 & x_0 & x_1 & x_2 & & \vdots\cr
\vdots && \ddots & \ddots & \ddots  \cr
0 & \ldots & & x_0 & x_1 & x_2 \cr}.
$$

\noindent{\sl Proof of \ref{rank 2 on P2}.\/}
Setting
$
\F=\coker \bigl(\O_{\P^2}(d-2)^{d-1}\rTo^\alpha \O_{\P^2}(d-1)^{d+1}\bigr)
$
we see that 
$\wedge^2 \F \iso \O_{\P^2}(3d-3)= \O_{\P^2}(3d)\otimes \omega_{\P^3}$.
Since $\F$ is a rank 2 vector bundle, 
$$
\F=\F^*\otimes \wedge^2\F 
= \F^*\otimes \O_{\P^2}(3d)\otimes \omega_{\P^3},
$$
so, {\it as a sheaf on the ambient space of the $d$-uple
embedding of $\P^2$,\/} $\F$ satisfies the duality hypothesis
of \ref{self-dual Ulrich}. Further, the given
presentation that $\F$ shows that $\F$ is $(d-1)$-regular
as a sheaf on $\P^2$,
and thus it is 0-regular on the ambient space
of the $e$-uple embedding for any $e\geq d-1$.
Thus \ref{self-dual Ulrich} shows that $\F$ is
an Ulrich sheaf on the $d$-uple embedding.\Box
  
The betti diagram of the Tate resolution of such a rank
2 sheaf $\F$ is given just after 
\ref{d-uple}. 
Instead of specifying $\alpha$, we could define
$\F$ by giving the $(d-1) \times 2(d-2)$ 
matrix $\beta$ of linear forms over $E$
that occurs at the end of the middle strand
of the Tate resolution. For the choice of $\alpha$ above
we get 
$$\beta=\pmatrix{
e_0 & e_1 & 0 & 0 && \ldots && 0 \cr
e_1 & e_2 & e_0 & e_1 & &  && \cr
0 & 0 & e_1 & e_2 & &&& \vdots\cr
\vdots &&&&\ddots  \ddots &&& 0 \cr
&&&&& &e_0 & e_1 \cr
0&&&\ldots&&0& e_1 & e_2 \cr },$$
and the vector bundle $\cE$ has a conic of maximal order jumping
lines. One can show by semi-continuity that $\beta$ can be
taken to be any sufficiently general 
$(d-1) \times 2(d-2)$ matrix of linear forms over $E$, but 
unlike for the matrix $\alpha$, we do not know how to recognize
when $\beta$ is sufficiently general to give rise to 
a Tate resolution of the right form.

\medskip \noindent {\bf Bundles on $\PP^3$} \medskip

\proposition{Instanton bundles} Suppose $d\geq 2$.
There exist rank 2 Ulrich sheaves
for the $d$-uple embedding of $\PP^3$ if and only if 
$d \not\equiv 0\  ({\rm mod}\ 3)$.

\proof
By Hartshorne and Hirschowitz
[1982] there exist  rank 2 vector bundles $\F$ with $c_1=0$ and
and natural cohomology on $\PP^3$ for any given $c_2$.
For 
$d \not\equiv 0\ ({\rm mod}\ 3)$  and
$c_2=(d^2-1)/3$ the sheaf $\F(d-2)$ is Ulrich for the d-uple embedding.
The converse follows from \ref{divisibility}. \Box

\remark{} The bundles $\F$ in the proof of the proposition
 are  called ``instanton bundles'', see Tikhomirov [1997],
because they satisfy the instanton conditions
$$ \F \hbox{ is stable of rank } 2,\, c_1(\F)=0 \hbox{ and } \H^1(\F(-2))=0.$$ 
Equivalently their linear monad $\LL(\F)$ has shape
$$ 0 \rTo \O(-1)^{c_2} \rTo \O^{2c_2+2} \rTo \O(1)^{c_2} \rTo 0.$$
Except for the 2-uple embedding,
it is an open problem us to find an explicit expression for
these rank 2 Ulrich sheaves. 

For the 2-uple embedding the rank 2 Ulrich sheaf is essentially unique:

\proposition{2-uple on P3} 
If $\E$ is the Null-correlation bundle on $\P^3$,
then $\F:=\E(-2)$ is,
up to automorphisms of $\P^3$, the unique
rank 2 Ulrich sheaf on the 2-uple embedding on $\P^3$,

\proof By \ref{d-uple}, $\F$ is an rank 2 Ulrich sheaf if and only if
the betti diagram of the Tate resolution of $\F$ has the form
$$\matrix{
***&64&35 & 16 & 5 &. &.&.&.&.&.&\dots\cr
\dots&.&.&.&.&1&.&.&.&.&.&\dots\cr
\dots&.&.&.&.&.&1&.&.&.&.&\dots\cr
\dots&.&.&.&.&.&.&5&16&35&64&*** \cr}$$  
with Hilbert polynomial $\chi \cF(t) = {1\over 3}(t+2)(t+4)(t+6)$.
(Here and henceforward, 
we replace each zero in a Betti diagram by a ``.'' to improve legibility.)
As proved in Okonek-Schneider-Spindler [1980],
the Null correlation bundle is determined (up to twist)
by its intermediate cohomology $\cF$
and the choice of a nondegenerate 2-form,
here given by the map in the middle of the Tate resolution.
Thus $\F$ must be a twist of the null correlation bundle,
the twist is determined by a comparison of Hilbert polynomials.
\Box

By \ref{divisibility} there is
no rank 2 bundle on $\P^3$ that is an Ulrich sheaf for the 
3-uple embedding. \ref{homogeneous Ulrich} gives
a homogeneous bundle of rank 9.
The following example gives a whole family
of rank 3 Ulrich bundles for this case.
These bundles give
determinantal B\'ezout formulas for the cube of the resultant
of 4 forms of degree 3 in 4 variables.

\example{3-uple on P3}
{\it A family of rank 3 vector bundles on $\P^3$ which
are Ulrich sheaves for the 3-uple embedding.}

By \ref{d-uple} $\F$ is an Ulrich sheaf
for the 3-uple embedding if and only if  the betti diagram of its
Tate resolution has the form
$$
\matrix{
\ldots&81&40&14 &  & . &. &.&.&.&.&.\cr
.&.&.&.&5&4&.&.&.&.&.&.\cr
.&.&.&.&.&.&4&5&.&.&.&.\cr
.&.&.&.&.&.&.&.&14&40&81&\ldots& \cr}
$$
Calculation shows that if we take a sufficiently general
$5\times 4$ matrix over the exterior algebra in 4 variables,
then its Tate resolution has this form.\Box

It follows at once from the definitions that
a sheaf on $\P^k$ becomes weakly Ulrich on the $d$-uple embedding
if and only if 
$$
\eqalign{
&h^0\F(-2d)=0; \cr
&h^i\F((-i-2)d)=0=h^i\F((-i+1)d)\quad 0<i<k-1;\ {\rm and} \cr
&h^k\F((-k+1)d)=0.
}
$$
{}From the form of the
cohomology diagram of the ``null correlation bundle'' 
on $\P^3$ given in the proof of \ref{2-uple on P3} we see that a twist of this bundle 
becomes weakly Ulrich on each $d$-uple embedding, and thus gives
a Pfaffian Stiefel formula for the 
of the resultant of 4 forms in 4 variables of any degree.
For any $d \ge 2$ the corresponding  2-term complex on 
$\GG(4,\H^0\cO_{\PP^3}(d))$ has the form
$$ 
0 \to \cO(-1)^b \oplus \cU^a \to 
\cO^b \oplus (\Lambda^3 \cU)^a \to 0
$$
with $a=d(d^2-4)/3$ and $b=2d(4d^2-4)/3$.

\medskip \noindent {\bf Bundles on $\PP^4$} \medskip

\example{Horrocks-Mumford} The Horrocks-Mumford bundle on $\PP^4$
has rank 2 and  Tate resolution
$$\matrix{ 
\ldots & 100 & 35 & 4 & . & . & . &.&.&.&.&.&. \cr
      . & . & 2 & 10 & 10 & 5 & . &.&.&.&.&.&. \cr
                  .&.&.&.&.&.&2 &.&.&.&.&.&. \cr
. &.&.&.&.&.&.&5&10&10&2&.&.\cr
. &.&.&.&.&.&.&.&.&4 &35 & 100 & \ldots \cr}$$
It gives rise to Pfaffian Stiefel formulas for $d=4,6,8$.

\example{d=2 Ulrich sheaf on P4} Suppose again that $k=4$,
and take $d=2$.
By \ref{divisibility}
any Ulrich sheaf
on the $2$-uple embedding of $\P^4$ has rank divisible by 8.
Consider a general map $E^3 \to E^5(-2)$. Its Tate resolution is
$$\matrix{ 
\ldots  &128 & 35 & . & . &.&.&.&. \cr
       . & . & . & 5 & . &.&.&.&. \cr
                  .&.&.&.&3 &.&.&.&. \cr
.&.&.&.&.&5&.&.&.\cr
.&.&.&.&.&.&35&128 &\ldots  \cr}$$
This gives a rank 8 Ulrich sheaf.
\medskip

\section{Surfaces} Surfaces

Throughout this section, $X$ denotes a nonsingular projective
surface over $K$, and we assume that $K$ has characteristic 0.
We study Ulrich sheaves on $X$.
We write $H$ for a hyperplane divisor and $K$ for a canonical
divisor on $X$.

In general it is rare to find an Ulrich line bundle on 
a surface; for example is easy to see that there are none on
the $d$-uple embedding of $\P^2$ when $d>1$. 
Thus we turn to rank 2 bundles.
By \ref{self-dual Ulrich},
If $\F$ is a rank 2 vector bundle on $X$ such that
$c_1 \F = 3H+K$ and $\F$ is 0-regular then $\F$ is
Ulrich by \ref{self-dual Ulrich}. We can obtain
a Pfaffian B\'ezout expression for
the Chow form from $\F$  on $X$. We will call such a
rank 2 bundle a {\it special rank 2 Ulrich bundle\/}. 

Many surfaces have no rank 2 Ulrich bundles. 
For example one can see by considering the
dimensions of the families that
the general surface $X$ of degree $d\ge 16$ in $\PP^3$ 
is not defined by the Pfaffian
of a $2d \times 2d$ skew symmetric linear matrix
(see Beauville [2000]).
Thus such a surface has no special rank 2 Ulrich bundle,
and because Pic $X=\ZZ$ for a general surface, every
rank 2 Ulrich sheaf would be special.

We are particularly interested in the case when $\F$ is
a blow-up of $\P^2$; then the Chow form is the resultant
of some ternary forms with some assigned base points
(that is, the vanishing of the Chow form determines when
these forms have an extra zero in common.)

\proposition{Ulrich on a surface}
\item{a)} Let $C$ be a smooth curve on $X$ of class $3H+K$ and let $\cL$ be a line
bundle on $C$ with
$$
\deg \cL = {1\over2} H.(5H+3K) + 2\chi\cO_X.
$$
If $\sigma_0,\sigma_1\in \H^0(\cL)$ define a base point free pencil
and  $\H^1 \cL(H+K)=0$
then the bundle $\F$ defined by the ``Mukai exact sequence''
$$
0 \rTo \F^* \rTo \cO^2_X \rTo^{(\sigma_0,\sigma_1)} \cL \rTo 0
$$ 
is a special rank 2 Ulrich bundle. 
\item{b)}Every special rank 2 Ulrich bundle on $X$
can be obtained from a Mukai sequence as in part a).

\proof a): We begin by proving that, under the hypotheses of part a),
the map
$$
(\H^0 \cO(H+K))^2 \rTo^{(\sigma_0,\sigma_1)} \H^0 \cL(H+K)\leqno{(*)}
$$
is an isomorphism. Using Riemann-Roch on $X$ and on $C$,
and the given degree of $\L$, we immediately compute
$\chi(\L(H+K))=2\chi(\O_X(H+K)$ and
$\chi(\L(2H+K))=2\chi(\O_X(2H+K)$. Our hypothesis that
$\L(H+K)$ is nonspecial implies that $\L(2H+K)$ is also
nonspecial. With this and the
Kodaira vanishing theorem on $X$, we
see that  $\chi$ is equal to $\h^0$ for
all four of these bundles. Thus it suffices to show that
the map $(*)$ is injective.

Since $C \sim 3H+K$
there is an exact sequence
$$
0\rTo \O_X(-2H)\rTo \O_X(H+K)\rTo O_C(H+K)\rTo 0,
$$
from which we see that the restriction map
$\H^0\cO_X(H+K)\cong \H^0\cO_C(H+K)$ is an injection.
By the base point free pencil trick there is a left exact sequence
$$
0\rTo \H^0 \cL^*(H+K) \rTo (\H^0 \cO_C(H+K))^2
\rTo^{(\sigma_0,\sigma_1)} \H^0\L(H+K),
$$
By the adjunction formula $K_C=(3H+2K)\mid_C$, so
our hypothesis and Serre duality
give $0=\h^1\L(H+K)=\h^0 \cL^*(2H+K)$,
whence $\h^0 \cL^*(H+K)=0$ as well. Thus
(*) is an injection.

We can now prove that $\F$ is Ulrich. 
The Mukai sequence implies that $\wedge^2\F = \O_X(3H+K)$,
so by
\ref{self-dual Ulrich} it suffices to show that
$\F$ is 0-regular.
Twisting the Mukai sequence by $H+K$ and using 
the preceding result together with Kodaira vanishing,
we see that $\H^1(\F^*(H+K))=\H^2(\F^*(H+K)) = 0$.
Serre duality now gives $\H^1(\F(-H))=0$ and $\H^0(\F(-H))=0$.
Since $\wedge^2\F = \O_X(3H+K)$,
we have $\F(-H)=\F^*(2H+K)$. By Serre duality
$\h^2(\F(-2H))=\h^0(\F^*(2H+K))=\h^0(\F(-H))=0$, 
and $\F$ is Ulrich as claimed.

b): Conversely, if $\F$ is a 
special Ulrich bundle of rank 2, then
two general sections $\tau_0,\tau_1$ of $\F$ become dependent
on a smooth curve $C$ of class $3H+K$. The cokernel of the 
induced map
$ 0 \to \F^* \to \oplus_1^2 \cO$
is a line bundle $\cL$ on $C$, generated by
global sections, so we obtain the Mukai sequence
$$ 
0 \to \F^* \to \oplus_1^2 \cO\rTo^{(\sigma_0,\sigma_1)} \L\to 0.
$$
By Serre duality, $\chi(\F^*(H+K))=\chi(\F(-H))$, which is
0 since $\F$ is Ulrich. Thus $\chi(\L(H+K))=2\chi(O_X(H+K))$.
Applying the Riemann-Roch theorems on $X$ and $C$
again, we obtain the
desired formula for the degree of $\L$.\Box

\corollary{del Pezzo} Suppose that the base field $k$ is
algebraically closed. If $X\subset \P^r$ is a del Pezzo surface, then
$X$ has a special rank 2 Ulrich Bundle. Thus there is a Pfaffian
B\'ezout formula for the resultant of 3 ternary cubics with
$d$ basepoints in general position.

\proof  In this case $K=-H$ and $C\sim 3H+K$ is a canonical curve
of genus $g=H^2+1=r+1$.
Any general line bundle of
degree 
$$
\deg \cL = {1\over2} H.(5H+3K) + 2\chi\cO_X=g+1.
$$ 
defines a nonspecial pencil. Thus we can apply
\ref{Ulrich on a surface} to get a special rank 2 Ulrich bundle on $X$.

The space of ternary cubics with $d$ general base points
has dimension $10-d$, so
it suffices to treat the case
of seven or fewer points. The linear series of
cubics with 7 assigned base points maps the plane two-to-one
onto itself, and the condition that three such cubics meet
in an extra point is the condition that three lines in the plane
meet in a point---a determinantal condition. 

For six or fewer assigned base points the resultant is
exactly the Chow form of the corresponding del Pezzo surface.
\Box

\corollary{} Let $C$ be a smooth curve of class $3H+K$ and let $\cL$ be a line
bundle on $C$ such that $\mid \cL \mid$ is a base point free pencil
of degree 
$
\deg \cL = {1\over2} H.(5H+3K) + 2\chi\cO_X.
$
The conditions of \ref{Ulrich on a surface} are satisfied
iff $\mid \cL \mid$ does not arise as a projection
from $\mid \cO_C(2H+K) \mid$.

\proof To say that $\mid \cL \mid$ arises as a projection
from $\mid \cO_C(2H+K) \mid$ means that 
of $\H^0(\cL^*(2H+K))\neq 0$. 
This space is Serre dual to $\H^1(\cL(H+K)=0$.
\Box

\remark{} Pencils which arise as projections correspond to
codimension 2 planes that are ${7 \over
 2}H.(H+K) + K^2 -2\chi\cO_X-$secant 
to $C \subset \PP H^0 \cO(2H+K)$. Every component of the variety
of such secants has dimension at least
$$ {1 \over 2}H.(H-K) + 4\chi \cO_X -4 -K^2, $$
and we might expect equality. On the other hand the variety of pencils
$\mid \cL \mid$ has dimension at least
$$
\rho(\cL) = 2 \deg \cL - g_C - 2 ={1 \over 2}H.(H-3K) + 4\chi \cO_X -1
-K^2.
$$
Thus we would
expect the existence of an $\cL$ which is not a projection,
and thus of a special rank 2 Ulrich bundle, in case 
$ 
H.K < 3.
$  
\medskip

\medskip
\noindent 
{\bf Resultants of ternary forms with base points} \medskip

\noindent
Consider $X=\PP^2(p_1,\ldots,p_e)$ the
blow up of the plane in $e$ distinct points and a very ample
divisor class $H=dL-\sum_{i=1}^e E_i$. Here $L$ denotes the class of a
line and the $E_i$ the exceptional divisors.
The Chow form of $X$ can be interpreted as the resultant of ternary
forms of degree $d$ with $e$ assigned base points.

\theorem{rational surfaces} Let the ground field be infinite. 
Let $E = \{p_1,\ldots,p_e\}$ be a collection of 
$e$ distinct points in $\PP^2$ and let $X = \PP^2(p_1,\ldots,p_e)$ 
be the blow up of $\PP^2$
in these points, embedded by the linear system $|dL-\sum E_i|$. If the 
homogenous ideal
$I_E$ of the points is generated in degree $d-1$ then $X$ has a special 
rank 2 Ulrich sheaf.

\proof Let $\eta : X \to \PP^2$ be the blow up. By \ref{Ulrich on a surface}
we have to construct a pencil $|\cL|$ of degree ${(d-1)(5d-4)\over 2} - e$ 
on a smooth curve of class $(3d-3)L-2\sum_iE_i$ on $X$ which 
satisfies $\H^1 \cL(H+K)=\H^1(\cL((d-3)L)=0$.
Let $C'=\eta(C) \subset \PP^2$ be the plane model. Every pencil on $C$ 
corresponds to a pencil of adjoint curves of degree a, say, with assigned
base points $F=q_1+\ldots+q_f$ on $C'$, that is
a pencil $\{\lambda A_0+\mu A_1\} \subset \H^0(\PP^2 ,\cI_{E\cup F}(a))$.
The pencil of plane curves might have additional base points
$G=r_1+\ldots+r_g$ away from $C'$. We have
$$ a^2= e + f +g. $$
In order that $|\cL|$ is not a projection from $|2H+K|$ we need $a>2d-3$.
We choose $a=2d-2$ so that we can deal with the fewest number of 
additional points $F$ and $G$.
To complete the construction we will choose $C'$ and $\L'=\eta_* \L$
simultaneously.

Take  $G=r_1+\ldots r_g$ as $g= { d \choose 2}$
general points in the plane disjoint from $E$. 
By the Hilbert-Burch theorem, see [Eisenbud, 1995, 20.4], $I_G$ is 
generated by the 
$d-1$ minors of a $d \times (d-1)$ matrix 
$\varphi_1: \cO_{\PP^2}(-1)^{d-1} \to \cO_{\PP^2}^d$ 
with linear entries, since $G$ impose independent conditions on 
forms of degree $d-2$. Since $I_E$ and $I_G$ are generated by forms of
degree $d-1$ the sheaf $\cI_{E\cup G}(2d-2))$ is globally generated.
Choose a general pencil
$$A_0,A_1 \in \H^0(\PP^2,\cI_{E\cup G} (2d-2)).$$
Then in our construction $F$ has to be the scheme defined by the ideal
$$I_F=(A_0,A_1):I_{E\cup G}$$
and will consist of $f$ simple points disjoint from $E \cup G$ by 
Bertini's theorem.
$C'$ and $\cL'$ are then presented as follows:
First note that the Hilbert-Burch matrix of
$I_{E\cup F}=(A_0,A_1):I_G$ is the $ d \times (d+1)$ matrix 
$\varphi_2: \cO_{\PP^2}(-1)^{d-1}\oplus \cO_{\PP^2}(-d+1)^2 \to \cO^d$
obtained from $\varphi_1$ by writing $A_0$ and $A_1$ as a linear 
combination of the 
generators of $I_G$, c.f [Peskine-Szpiro, 1974]. Since
$I_{2E \cup F} \subset I_{E \cup F}$ we can obtain the equation of $C'$ as 
 a determinant of a
matrix $\varphi_3 : \cO_{\PP^2}(-1)^{d-1}\oplus \cO_{\PP^2}(-d)^2 \to \cO_{\PP^2}^{d+1}$ 
obtained from $\varphi_2$
by adding a column. The transposed matrix twisted gives $\cL'$:
$$
0 \to \cO_{\PP^2}(-d+1)^{d+1} \rTo^{\varphi_3^t} \cO_{\PP^2}(-d+2)^{d-1}\oplus 
\cO_{\PP^2}^2 \rTo  \cL' \to 0.
$$
For a general choice $C' \in |(3d-3)L-2E-F|$ the curve $C'$ will have only 
ordinary double points in $E$. 
For example we could simply take all entries
of degree $d-1$ in the matrix $\varphi_3$ general elements in 
$(I_E)_{d-1}$ and 
all the linear forms general.
 
Locally around a point $p_i$ of $E$ the sheaf
$\L'$ defined by the sequence above has a stalk $\cL_{p_i}$ which is minimally
generated by two elements, since $A_0$ and $A_1$ intersect 
transversally at $p_i$.
So $\cL'=\eta_* \cL$ for some line bundle $\cL$ on $C$.
   
Since the additional generators of $\sum_m \H^0(\eta_* \cL(m))$ 
are of degree $d-2$
and $\H^0(X,\cO_X(H+K) \cong \H^0(\PP^2, \cO(d-3)$ we see that 
the desired isomorphism $(*)$
$\H^0(\cO_X(H+K))^2 \cong \H^0(\cL(H+K))$ from the proof of 
\ref{Ulrich on a surface} holds. This completes the proof.\Box

\corollary{} There exists a Pfaffian B\'ezout formula for ternary forms 
of degree $d$ with $e$ assigned base points if the ideal of the points
is generated in degree $d-1$. \Box

\remark{larger range} Our computations suggest that
the construction of the rank 2
Ulrich sheaf above, and hence the construction of a B\'ezout formula
for forms with base points,
works for a set of points $E$ even under the weaker hypothesis
that
 $I_E$ is generated in degree $d$. For example, if $E$ consists
of $e\le {d+2 \choose 2} -6$ general points, there 
should be plenty of room to arrive at a nodal $C'$ in the construction.
 
\section{Appendix} Appendix by Jerzy Weyman: $\Ext(\wedge^qU, \wedge^pU)$

In this appendix we will prove part b of \ref{Hom computation}, and 
also prove a complementary statement about the higher cohomology. In 
characteristic 0 these statements follow from the Bott vanishing
theorem, but we prove them over a field of arbitrary characteristic:

\theorem{all char hom} Let $\GG_l$ be the Grassmannian of codimension
$l$ planes in a vector space $W$ with dual $V=W^*$ over a field $K$
of arbitrary characteristic, and let $U$ be the
tautological
$l$-sub bundle of $W\times \GG_l$. For $0\leq p,q,\leq l$ we have
$$
\Hom(\wedge^qU, \wedge^pU)=
\cases{0;&{\rm if}
$p>q$\cr\bigwedge^{q-p}V ;&\rm{otherwise.}}
$$
Moreover $\Ext^i(\wedge^qU, \wedge^pU)=0$ for $i>0$ and all $p,q$.

Let $GL=GL_K(W)$ be the
general linear group. We write $Q$ for the tautological quotient
bundle $Q=W/U$ on $\GG_l$.
If $\lambda=(\lambda_1,\dots,\lambda_v)$ is a nondecreasing sequence of positive
integers (a highest weight for $GL$) then we write $L_\lambda W$ for the Schur
module corresponding to the highest weight $\lambda$. 
We may extend this notation to any nondecreasing sequence of
integers $\lambda$ (dominant integral weight)
using the formula
$L_\lambda W=L_{\mu^\prime} W \otimes (\bigwedge^v W)^{\otimes \lambda_v}$
where 
$\mu^\prime$ is the partition conjugate to $\mu = (\lambda_1
-\lambda_v ,\ldots ,\lambda_{v-1}-\lambda_v ,0)$.
The proof of \ref{all char hom} rests on the following
facts:

\lemma{A1} The tensor product  
$\bigwedge^p {U}\otimes\bigwedge^q {U}^* $ has a filtration with the
associated graded object
$$\oplus_{a+b=p-q, 0\le a\le p,0\le b\le q, a+b\le l} L_{(1^a , 0^{l-a-b},
(-1)^b )}{U}.$$

\lemma {A2} a) If if $a>0$ then all cohomology groups of the vector bundles 
$L_{(1^a , 0^{l-a-b},(-1)^b )}{U}$ 
are zero.\hfill\break
b) All higher cohomology groups of the bundle 
$L_{( 0^{l-b}, (-1)^b)}{U}$ 
are zero and
$$
H^0 (\GG_l ,L_{( 0^{l-b}, (-1)^b )}{U} )=\bigwedge^b W^*.
$$

\noindent{\sl Proof of \ref{A1}} This is a standard fact on good
filtrations (see Donkin [1985]) that the tensor product of Schur
modules has the filtration with associated graded being a direct sum of
Schur modules. The multiplicities of the Schur modules
occurring are the same as in characteristic zero, and we can get the result
by Littlewood-Richardson rule, using the
isomorphism $\bigwedge^q {U}^* =\bigwedge^{l-q}{U}\otimes
\bigwedge^l {U}^*$.\Box

\noindent{\sl Proof of \ref{A2}} Let $\lambda = (\lambda_1 ,\ldots
,\lambda_v )$ be an $v$-tuple of integers. Consider the full
flag variety and the tautological subbundles ${U}_i$ of rank $i$ on it.
We denote by ${\cal L}(\lambda )=
\otimes_{1\le i\le v} ({U}_i /{U}_{i-1})^{-\lambda_i}$ 
the line bundle on the full
flag variety ${GL}/{B}$, where $B$ is the Borel subgroup. Then we have

\lemma {A3} a) If $\lambda$ is a dominant integral weight, then the higher
cohomology groups of ${\cal L}(\lambda )$ vanish
and
$$
H^0 ({GL}/{B} ,{\cal L}(\lambda ))= L_\lambda W.
$$
b) Let us assume that for some $i$ we have $\lambda_i=\lambda_{i-1}
+1$. Then all cohomology groups of ${\cal
L}(\lambda )$ vanish.

Now part b) of \ref{A2} follows from part a) of \ref{A3}. To prove part a) of
\ref{A2} we consider the natural projection 
$\eta: {GL}/{B}\rightarrow \GG_l$. 
We
observe that by Kempf's Vanishing Theorem (see Jantzen [1987]) in the relative
setting we have
$L_{(1^a , 0^{l-a-b}, (-1)^b )}{U} =\eta_* ({\cal L}(0^{v-l},1^a ,
0^{l-a-b}, (-1)^b )) $  with higher direct images
${\cal R}^i\eta_* ({\cal L}(0^{v-l},1^a , 0^{l-a-b}, (-1)^b )) $ being zero
for $i>0$. Since by lemma 3 b) we know that all
cohomology groups of ${\cal L}(0^{v-l},1^a , 0^{l-a-b}, (-1)^b )$ are zero,
by the spectral sequence of the composition we are
done.\Box

\noindent{\sl Proof of \ref{A3}} The part a) is just Kempf's Vanishing
Theorem. Part b) follows from the following
consideration. Let $P(i)$ be a parabolic subgroup such that the
corresponding homogeneous space is a flag variety of
flags of dimensions $(1,2,\ldots ,i-1,i+1,\ldots ,v-1,v)$. The projection
$\rho :{GL}/{ B}\rightarrow {GL}/P(i)$
allows to identify ${GL}/{B}$ with the projectivization 
$\PP({U}_{i+1}/{U}_{i-1})$. The bundle 
${\cal L}(\lambda )$ is of the form 
$\rho^* ({\cal M})\otimes {\cal O}_{\PP({U}_{i+1}/{U}_{i-1})}(-1)$
because all the factors in the definition of ${\cal L}(\lambda_1 ,\ldots
,\lambda_v )$ except of the $i$-th and $i+1$-st are
induced from ${GL}/P(i)$. Therefore by the Serre's Theorem (in
relative setting) and by the projection formula we see that
all higher direct images ${\cal R}^i\rho_* ({\cal L}(\lambda_1 ,\ldots
,\lambda_v ))$ are zero. This implies part b). \Box

\references
\frenchspacing
\parindent=0pt

\medskip
\noindent
C.~d'Andrea and A.~Dickenstein: Explicit formulas for the multivariate
resultant.  Preprint, 2000.

\medskip
\noindent
B.~Ang\'eniol and M.~Lejeune-Jalabert: {\sl Calcul diff\'erentiel et
classes  caract\'eristiques en g\'eom\'etrie alg\'ebrique.\/}
Travaux en Cours 38, Hermann, Paris, 1989. 

\medskip
\noindent
A.~Beauville: Determinantal hypersurfaces. Dedicated
to William Fulton on the occasion of his 60th birthday.
Mich. Math. J. 48 (2000) 39--64.

\medskip
\noindent
A.~Beilinson: Coherent sheaves on $\P^n$ and problems
of linear algebra. 
Funct. Anal. and its Appl. 12 (1978) 214--216.
(Trans. from Funkz. Anal. i. Ego Priloz 12 (1978) 68--69.)

\medskip
\noindent
I.N.~Bernstein, I.M.~Gel'fand and S.I.~Gel'fand:
Algebraic bundles on $\P^n$ and problems of linear algebra.
Funct. Anal. and its Appl. 12 (1978) 212--214.
(Trans. from Funkz. Anal. i. Ego Priloz 12 (1978) 66--67.)

\medskip
\noindent
E.~B\'ezout: {\sl Th\'eorie g\'en\'erale des \'equation alg\'ebraiques,}
Pierres, Paris 1779. 

\medskip
\noindent
J.~Brennan, J.~Herzog, and B.~Ulrich: 
Maximally generated Cohen-Macaulay modules.
Math.~Scand.~61 (1987) 181--203.
  
\medskip
\noindent
R.-O. Buchweitz, D. Eisenbud, and J. Herzog:
Cohen-Macaulay modules on quadrics, by  In {\sl Singularities,
representation of algebras, and vector bundles (Lambrecht, 1985),\/}
Springer-Verlag Lecture Notes in Math. 1273 (1987) 96--116.

\medskip
\noindent
R.-O. Buchweitz: Maximal Cohen-Macaulay modules and Tate-cohomolgy
over Gorenstein rings, unpublished manuscript, 1987.

\medskip
\noindent
R.-O.~Buchweitz and F.-O.~Schreyer: Complete intersections of two quadrics, 
hyperelliptic curves and their Clifford Algebras, manuscript, 2002.

\medskip
\noindent
A.~Cayley: On the theory of elimination. Cambridge and Dublin
Mathematical J.~ 3 (1848) 116--120.

\medskip
\noindent
S.~Donkin: {\it Rational representations of algebraic groups.
Tensor products and filtration. \/}
Lecture Notes in Math.~1140.
Springer-Verlag, Berlin, 1985. 

\medskip
\noindent
D.~Eisenbud: {\sl Commutative Algebra with a View Toward 
Algebraic Geometry.} Springer Verlag, 1995.

\medskip
\noindent
D. Eisenbud: Linear sections of determinantal varieties.
American J. Math. 110 (1988) 541--575.

\medskip
\noindent
D.~Eisenbud and S.~Goto: Linear free resolutions and minimal
multiplicity.  J. Alg. 88 (1984) 89--133.

\medskip
\noindent
D. Eisenbud, G.~Fl\o ystad and F.-O. Schreyer: Sheaf cohomology and free
resolutions over the exterior algebra. Preprint (2000).

\medskip
\noindent
W.~Fulton: Young tableaux : with applications to representation theory and
geometry. 
London Mathematical Society student texts 35,
Cambridge University Press 1997.

\medskip
\noindent
I.~M.~Gelfand, M.~Kapranov, and A.~Zelevinsky: 
{\sl Discriminants, resultants, and multidimensional determinants\/}
Birkh\"auser, Boston, 1994.

\medskip
\noindent
D.~Grayson and M.~Stillman: Macaulay2.\hfill\break
{\tt http://www.math.uiuc.edu/Macaulay2/} (1993--).

\medskip
\noindent
D.~Hanes: Special condition on maximal Cohen-Macaulay modules, and 
applications to the theory of multiplicities. Thesis, Michigan, 2000.

\medskip
\noindent
R. Hartshorne: {\sl Algebraic Geometry.} Springer Verlag, 1977.

\medskip
\noindent
R.~Hartshorne and A.~Hirschowitz: Cohomology of the general instanton bundle. 
Ann. scient. \'Ec. Norm. Sup. a s\'erie 15 (1982) 365--390.

\medskip
\noindent
J.~Herzog, B.~Ulrich, and J.~Backelin: Linear maximal Cohen-Macaulay
modules over strict complete intersections.
J.~Pure and Appl.~Alg.~71 (1991) 187--202.

\medskip
\noindent
G. Horrocks, D. Mumford: A rank $2$ vector bundle on $\P^4$ with $15,000$
symmetries. Topology 12 (1973) 63--81.

\medskip
\noindent
B.~Iversen: 
{\it Cohomology of Sheaves.\/}
Springer-Verlag, New York, 1986.

\medskip
\noindent
J.~C.~Jantzen:
{\it Representations of algebraic groups.\/}
Pure and Applied Mathematics, 131, 
Academic Press, Boston, MA, 1987. 

\medskip
\noindent
J.-P.~Jouanolou: Aspects invariants de l'Žlimination.
Adv. Math. 114 (1995) 1--174.

\medskip
\noindent
M.~Kline: {\it Mathematical thought from ancient to modern times}
Oxford University Press, New York 1972.

\medskip
\noindent
F.~Knudsen and D.~Mumford: The projectivity of the moduli
space of stable curves I: Preliminaries on ``det'' and ``div''.
Math.~Scand.~39 (1976) 19--55.

\medskip
\noindent
G.W.~Leibniz: {\sl Mathematische Schriften} herausgegeben von C.I.~Gerhardt,
Letter to l`Hospital 28 April 1693, volume II, pp. 239; Vorwort volume VII, pp. 5.
Georg Olms Verlag Hildesheim, New York 1971.

\medskip
\noindent
Ch. Okonek, M. Schneider, and H. Spindler: {\sl Vector Bundles on
Complex Projective Spaces.\/} Birkh\"auser, Boston 1980.

\medskip
\noindent
C.~Peskine and L.~Szpiro: Liaison des vari\'eti\'es alg\'ebriques I. 
Invent. Math. 26 (1974) 271--302

\medskip
\noindent
R.~Stanley: Theory and application of plane partitions. Studies in Appl. Math.  
1 (1971) 167--87 and 259--279.

\medskip
\noindent
R.~G.~Swan: $K$-theory of quadric hypersurfaces. Ann. of Math.
122 (1985) 113--153.

\medskip
\noindent
J.J.~Sylvester: A method of determining by mere inspection the derivatives from two 
equations of any degree. Philosophical Magazine XVI (1840) 132--135.
Memoir on the dialytic method of elimination. Part I. 
Philsophical Magazine XXI (1842) 534--539. Reprinted in:
{\sl The collected Mathematical Papers of James Joseph Sylvester, Volume I.}
Chelsea New York 1973, reprint of the Cambridge 1904 edition.

\medskip
\noindent
A.~A.~Tikhomirov: The main component of the moduli space of 
mathematical instanton vector bundles on $\PP^3$.
Journal of the Mathematical Sciences 86 (1997) 3004--3087.
(Translation from the Russian of {\it Itogi Nauki i Tkhniki.
Seriya Sovremennaya Matematika i Ee Prilozheniya. Tematischeskie Obzory}
Vol 25, Algebraic Geometry-2, 1995.)

\medskip
\noindent
B.~Ulrich: Gorenstein rings and modules with high numbers of generators.
Math.~Z.~188 (1984) 23--32.

\medskip
\noindent
J.~Weyman and A.~Zelevinsky: Determinantal formulas for multigraded resultants.
J.~Alg.~Geom.~3 (1994) 569--598.

\bigskip
\vbox{\noindent Author Addresses:
\smallskip
\noindent{David Eisenbud}\par
\noindent{Department of Mathematics, University of California, Berkeley,
Berkeley CA 94720}\par
\noindent{eisenbud@math.berkeley.edu}
\smallskip
\noindent{Frank-Olaf Schreyer}\par
\noindent{FB Mathematik, Universit\"at Bayreuth
D-95440 Bayreuth, Germany\par}
\noindent{schreyer@btm8x5.mat.uni-bayreuth.de}\par
\smallskip
\noindent{Jerzy Weyman}\par
\noindent{Department of Mathematics, Northeastern University, Boston MA 02115\par}
\noindent{weyman@neu.edu}\par
}
\rightline {Revised November 4, 2001}

\end